\documentclass[11pt,leqno]{amsart}
\usepackage{amssymb,amsmath,amsthm,amsfonts,color}
\usepackage{enumerate}
\usepackage{graphicx}

\usepackage[dvipsnames]{xcolor}

 \DeclareMathOperator{\Imag}{Im}
\DeclareMathOperator{\spa}{span} 
 
\DeclareMathOperator{\length}{length}
\DeclareMathOperator{\diam}{diam}
\DeclareMathOperator{\dist}{dist}\DeclareMathOperator{\id}{id}
\DeclareMathOperator{\Modd}{Mod}\DeclareMathOperator{\arccosh}{arccosh}
\DeclareMathOperator{\essinf}{essinf}
\DeclareMathOperator{\Wh}{Wh}
\DeclareMathOperator{\Bd}{Bd}

\newcommand{\field}[1]{\mathbb{#1}}
\newcommand{\N}{\field{N}}                      
\newcommand{\R}{\field{R}}                      
\newcommand{\Sph}{\field{S}}                    
\newcommand{\C}{\field{C}}                      
\newcommand{\Heis}{\field{H}}                    

\newcommand{\eps}{\epsilon}

\newcommand{\loc}{{\scriptstyle{loc}}}

\newcommand{\dR}{{\scriptstyle{dR}}}

\newcommand{\cB}{{\mathcal B}}

\newcommand{\cH}{{\mathcal H}}

\newcommand{\cL}{{\mathcal L}}

\newcommand{\cW}{{\mathcal W}}

\newcommand{\bp}{{\mathbf p}}
\newcommand{\bq}{{\mathbf q}}
\newcommand{\bA}{{\mathbf A}}

\newcommand{\norm}[1]{\lVert #1 \rVert}

\def\Barint_#1{\mathchoice
          {\mathop{\vrule width 6pt height 3 pt depth -2.5pt
                  \kern -8pt \intop}\nolimits_{#1}}%
          {\mathop{\vrule width 5pt height 3 pt depth -2.6pt
                  \kern -6pt \intop}\nolimits_{#1}}%
          {\mathop{\vrule width 5pt height 3 pt depth -2.6pt
                  \kern -6pt \intop}\nolimits_{#1}}%
          {\mathop{\vrule width 5pt height 3 pt depth -2.6pt
                  \kern -6pt \intop}\nolimits_{#1}}}

\theoremstyle{plain}
\newtheorem{theorem}{Theorem}

\newtheorem{conjecture}[theorem]{Conjecture}

\theoremstyle{definition}

\newtheorem{question}{Question}

\numberwithin{theorem}{section} 
\numberwithin{equation}{section}

\usepackage[hidelinks]{hyperref}
\hypersetup{ 
  colorlinks   = true,            
  urlcolor     = {blue!50!black}, 
  linkcolor    = {blue!50!black}, 
  citecolor   = {red!50!black}    
}

\title{The Heinonen--Semmes problems after thirty years}
\date{\today}

\author{Guy C. David}
\address{Department of Mathematical Sciences, Ball State University, Muncie, IN, USA}
\email{gcdavid@bsu.edu}

\author{Pekka Pankka}
\address{Department of Mathemetics and Statistics, University of Helsinki, Finland }
\email{pekka.pankka@helsinki.fi}

\author{Jeremy T. Tyson}
\address{Department of Mathematics, University of Illinois at Urbana-Champaign, Urbana, Illinois, USA}
\email{tyson@illinois.edu}

\begin{document}
\maketitle

\begin{abstract}
We survey the current status of the questions posed by Juha Heinonen and Stephen Semmes in `Thirty-three yes or no questions about mappings, measures, and metrics' ({\it Conf.\ Geom.\ Dynam.}, 1997).
\end{abstract}


{
    \hypersetup{linkcolor=black}
    \setcounter{tocdepth}{1}
    \tableofcontents
}

\section{Introduction}\label{sec:intro}

In 1997, Juha Heinonen and Stephen Semmes published an article \cite{hs:questions} in the first volume of the AMS journal {\it Conformal Geometry and Dynamics}\footnote{In fact, \cite{hs:questions} was the very first article published by {\it Conformal Geometry and Dynamics}.} entitled `Thirty-three yes or no questions about mappings, measures, and metrics'. As noted by the authors, the problems in this list were mostly fairly new at the time of compilation. The problems reside primarily within the fields of geometric measure theory, metric geometry, geometric topology and quasiconformal geometry, and often cross boundaries between these fields. Nowadays, the majority of these problems would be categorized under the heading of {\it analysis in metric spaces}. Analysis in metric spaces was only just beginning to emerge as an independent research area around the time when this article was published, and indeed this problem list was one of the first influential papers which helped to direct the future evolution of the subject. 

The purpose of this article is to survey the current status of the thirty-three Heinonen--Semmes problems, and also to offer brief comments about the role which this list of problems played in the development of analysis in metric spaces over the past thirty years. 

Our focus on the Heinonen--Semmes problem list is not intended to imply that these are the only interesting problems in the field. The subject of analysis in metric spaces has grown dramatically over the intervening decades, and now involves novel interactions with other research areas which were essentially unanticipated in the late 1990s. Also, there have been other important developments in the field not inspired by or related to this list of problems. Nevertheless, we believe it is beneficial -- at this moment in the history of the subject -- to take a historical look at the role played by this problem list.

In the introduction of \cite{hs:questions}, the authors write: ``None of [these problems] seem easy and some are likely to be very difficult. The formulation of each problem is such that it can be answered by one word only: either {\em yes} or {\em no} \dots We offer no conjectures or guesses.'' The decision not to outline potential approaches to any of these problems, nor to make any predictions as to which way the answer would go, could be regarded as unusual for a problem list of this type. However, in our opinion, this strategy was well-chosen. This choice made it easier for future generations of researchers to approach these problems {\it ab novo} and with a fresh perspective. In fact, as we observe below, the eventual solutions to several of the problems went against the general expectation at the time.

For the reader's convenience, we follow the original ordering of the problems in \cite{hs:questions}. We have grouped the problems into several broadly defined categories, while acknowledging that the interdisciplinary character of many problems makes our choice somewhat arbitrary. Within each category, we first state a group of problems verbatim from the original problem list, along with the current status of the problem. We then recall basic definitions and concepts needed to understand the problem statements. Finally, we provide more detail about the actual solution of each problem, including a precise statement of the relevant theorem(s). In cases where the problem remains unsolved (to the best of our knowledge), we indicate when appropriate any extant partial results. 

We conclude this survey with a brief overview of various directions within the broad heading of analysis in metric spaces, focusing particularly on topics which have arisen since the Heinonen--Semmes problem list was published. We emphasize in particular how those problems, and efforts undertaken towards their solutions, led to the development of connections between this field and related areas.

At the time of preparation of this survey, and to the best of our knowledge, $18$ questions have been solved (possibly not in all dimensions), with $15$ questions remaining open.

\subsection*{Acknowledgments}
Numerous individuals have contributed to the development of analysis in metric spaces over the past thirty years, and we benefited from the input and advice of many colleagues in the preparation of this manuscript. In particular, we are deeply grateful to those colleagues who joined us for an informative roundtable discussion during the June 2025 {\it Quasiweekend} conference held in Helsinki, some of whom also provided us with feedback on an early draft of this paper. We acknowledge in particular Christopher Bishop, Mario Bonk, Sylvester Eriksson-Bique, Hrant Hakobyan, Dimitrios Ntalampekos, Kai Rajala, Matthew Romney, Raanan Schul, Stefan Wenger, and Jang-Mei Wu for their comments. 

The quote in section 5.2, taken from the review \cite{hei:review} written by Juha Heinonen and published in Math Reviews as MR1731465, is included by permission from the American Mathematical Society.

\smallskip

GCD acknowledges support from the National Science Foundation under grant no.\ DMS-2054004. PP was partially supported by Research Council of Finland Center of Excellence FiRST. JTT acknowledges support from the Simons Foundation under Simons Collaboration Grant 852888.

\subsection*{Notation and conventions.} We review here basic notation and terminology. Concepts and ideas relevant for individual questions will be recalled in the section where those questions are discussed.

A {\it metric measure space} is a triple $(X,d,\mu)$, where $(X,d)$ is a metric space and $\mu$ is a Borel regular measure on $(X,d)$ which assigns a positive and finite value to each metric ball $B(x,r)$. We denote by $B(x,r)$ the open ball with center $x$ and radius $r>0$. The measure $\mu$ is {\it doubling} if there exists a constant $C>0$ so that $\mu(B(x,2r)) \le C \mu(B(x,r))$ for every metric ball $B(x,r)$ in $(X,d)$. A metric measure space $(X,d,\mu)$ is said to be {\it Ahlfors $Q$-regular}, $Q>0$, if there exists a constant $C>0$ so that $C^{-1} r^Q \le \mu(B(x,r)) \le C r^Q$ whenever $x \in X$ and $0<r\leq\diam(X)$. It is known that if $(X,d)$ admits a $Q$-regular measure $\mu$, then the $Q$-dimensional Hausdorff measure $\mathcal{H}^Q$ is also $Q$-regular and $\mu$ is comparable to $\mathcal{H}^Q$. Every $Q$-regular metric measure space with $Q>0$ is doubling. In the case when $X = \R^n$, we write $\cL^n$ for the Lebesgue measure and denote the element of integration by $dx = d\cL^n(x)$.

$L^p(X) = L^p(X,\mu)$ denotes the space of $p$-integrable functions on $X$, while $L^p_\loc(X)$ denotes the space of {\it locally $p$-integrable functions} on $X$. That is, $u \in L^p_\loc(X)$ if $u \in L^p(B(x,r))$ for every ball $B(x,r)$ in $X$. 

A mapping $f\colon (X,d) \rightarrow (Y,d')$ between metric spaces is \emph{Lipschitz} with constant $L>0$ if 
$$ 
d'(f(x),f(y)) \leq Ld(x,y) \qquad \mbox{for all $x,y \in X$,}
$$
and is \emph{bi-Lipschitz} with constant $L\ge 1$ if
\begin{equation}\label{eq:BL}
L^{-1}d(x,y) \leq d'(f(x),f(y)) \leq Ld(x,y) \qquad \mbox{for all $x,y \in X$.}
\end{equation}
Every bi-Lipschitz mapping is necessarily an embedding.

For a subset $A \subset (X,d)$, the {\it diameter} of $A$ is $\diam(A) = \sup\{d(x,y):x,y \in A\}$. A continuous mapping $\gamma:I \to X$ from an interval $I$ is a {\it curve}; we often identify the map $\gamma$ with its image $\gamma(I) \subset X$. A curve $\gamma$ defined on a compact interval $I$ is {\it rectifiable} if it has finite length, where $\length(\gamma)$ is defined to be the supremum of the values $\sum_{j=1}^m d(\gamma(t_j),\gamma(t_{j-1}))$ taken over all partitions $a = t_0 < t_1 < \cdots < t_{m-1} < t_m = b]$ of the interval $I=[a,b]$. A metric space $(X,d)$ is {\it geodesic} if any two points $x,y \in X$ can be joined by a rectifiable curve $\gamma$ such that $\length(\gamma) = d(x,y)$.

\section{Quasiconformal uniformization and parametrization}\label{sec:qc-uniformization-parametrization}

\begin{question}\label{q:1}
\it Is every strong $A_\infty$ weight in $\R^2$ comparable to the Jacobian of a quasiconformal map $f:\R^2 \to \R^2$?
\end{question}

The answer is NO, by work of Laakso \cite{laa:a-infinity}.

\begin{question}\label{q:2}
\it Is every strong $A_1$ weight in $\R^n$, $n \ge 2$, comparable to the Jacobian of a quasiconformal map $f:\R^n \to \R^n$?
\end{question}

The answer is NO, by a result of of Bishop \cite{bis:a-1}.

\begin{question}\label{q:3}
\it If $(\Sph^2,d)$ is both linearly locally contractible and Ahlfors $2$-regular, is it then quasisymmetrically equivalent to $\Sph^2$?
\end{question}

The answer is YES. This was first shown by Bonk and Kleiner in \cite{bk:spheres}. More recently, there have been a number of alternative approaches: by K. Rajala \cite{raj:spheres}, by Lytchak and Wenger \cite{lw:spheres}, by Meier and Wenger \cite{meierwenger}, and by Ntalampekos and Romney \cite{nr:poly}.

\begin{question}\label{q:4}
\it If $(\Sph^2,d)$ is both linearly locally contractible and Ahlfors $2$-regular, is it then bi-Lipschitz equivalent to $\Sph^2$?
\end{question}

The answer is NO. This follows from the ``No'' answer to Question \ref{q:1} given by Laakso in \cite{laa:a-infinity}. Laakso's example does not even admit a bi-Lipschitz embedding into any uniformly convex Banach space.

\begin{question}\label{q:5}
\it If $(\Sph^n,d)$, $n \ge 2$, is linearly locally contractible, Ahlfors $n$-regular, and quasisymmetrically three point homogeneous, is it then quasisymmetrically equivalent to $\Sph^n$?
\end{question}

The answer is YES for $n=2$ by a result of Bonk and Kleiner in \cite{bk:rigidity}, but the question remains open in dimensions $n\ge 3$.

\begin{question}\label{q:6}
\it Do smooth $n$-dimensional chord arc surfaces with small constant in $\R^{n+1}$, $n \ge 3$, admit quasisymmetric parametrizations?
\end{question}

This question remains open.

\subsection{Definitions and notation}\label{subsec:qc-uniformization-parametrization-definitions-and-notation}

Assume that $(X,d,\mu)$ is an Ahlfors $Q$-regular metric measure space. A nonnegative function $w \in L^1_\loc(X,\mu)$ is called a {\it weight} on $X$. A weight $w$ said to be a {\it strong $A_\infty$ weight} if the associated measure $d\nu := w d\mu$ is a doubling measure on $(X,d)$, and the induced quasi-distance function $D_w(x,y) := (\int_{B_{x,y}} w \, d\mu )^{1/Q}$ is comparable to a metric, i.e., there exists a constant $C>0$ and a metric $d'_w$ on $X$ so that
\begin{equation}\label{eq:dwprime}
C^{-1} d'_w(x,y) \le D_w(x,y) \le C d'_w(x,y) \qquad \text{for all } \, x,y \in X.
\end{equation}
Here and henceforth, for $x,y \in X$, we set
\begin{equation}\label{eq:B-x-y}
B_{x,y} := B(x,d(x,y))\cup B(y,d(x,y)).
\end{equation}
In case $w$ is a strong $A_\infty$ weight, the measure $d\nu := w d\mu$ is also called a {\it metric doubling measure}. 

A weight $w$ on $(X,d,\mu)$ is an {\it $A_1$ weight} if there exists $C>0$ so that 
$$
r^{-Q} \int_{B(x,r)} w \, d\mu \le C \essinf_{B(x,r)} w
$$
for all metric balls $B(x,r)$ in $X$. Here $\essinf_E g$ denotes the essential infimum of a function $g$ on a measurable set $E$ in $X$.

A homeomorphism $f:\Omega \to \Omega'$ between domains in $\R^n$, $n \ge 2$, is said to be {\it $K$-quasiconformal} for some $K \ge 1$ if $f$ lies in the local Sobolev space $W^{1,n}_\loc(\Omega:\R^n)$ and the weak differential of $f$ satisfies the pointwise inequality $||Df(x)||^n \le K \det Df(x)$ a.e.\ in $\Omega$. Here $Df(x)$ denotes the matrix of partial derivatives of the components of $f$, $||\bA||$ denotes the operator norm of an $n\times n$ matrix $\bA$, and $\det Df(x)$ denotes the Jacobian determinant of $f$. Alternatively, one may define a {\it (metrically) $H$-quasiconformal map} for some $H \ge 1$ to be an orientation preserving homeomorphism $f:\Omega \to \Omega'$ for which $H_f(x) \le H$ for all $x \in \Omega$, where $H_f(x)$ denotes the {\it linear dilatation} of $f$ at $x$, given by $H_f(x) = \limsup_{r \to 0} \max \{ |f(x)-f(y)|:|x-y|=r \} / \min \{ |f(x)-f(z)|:|x-z|=r \}$. It is important to note that Sobolev membership is not assumed in the metric definition of quasiconformality; part of the content of the statement that the two definitions are equivalent lies in the fact that a homeomorphism which satisfies the metric definition of quasiconformality necessarily lies in the local Sobolev space $W^{1,n}_\loc$ and is weakly differentiable a.e.

A homeomorphism $f:(X,d) \to (X',d')$ between metric spaces is said to be {\it $\eta$-quasisymmetric}, for some increasing homeomorphism $\eta:[0,\infty) \to [0,\infty)$, if for any $x,y,z \in X$ and $t>0$,
\begin{equation}\label{eq:qs}
d(x,y) \le t d(x,z) \qquad \text{if and only if} 
\qquad d'(f(x),f(y)) \le \eta(t) d'(f(x),f(z)).
\end{equation}
Every orientation preseving $\eta$-quasisymmetric homeomorphism $f$ between domains in $\R^n$, $n \ge 2$, is $K$-quasiconformal with $K = \eta(1)$, and every $K$-quasiconformal homeomorphism is locally $\eta$-quasisymmetric for some distortion function $\eta$ which depends only on $K$ and $n$. Moreover, a homeomorphism $f:\R^n \to \R^n$, $n \ge 2$ is quasiconformal if and only if it is quasisymmetric. It is easy to see that every bi-Lipschitz mapping (see \eqref{eq:BL} for the definition) is necessarily a quasisymmetry.

A metric space $(X,d)$ is said to be {\it linearly locally contractible} if there exists $C>0$ so that each metric ball $B(x,r)$ with $r<1/C$ can be contracted to a point inside the ball $B(x,Cr)$. The space $(X,d)$ is said to be {\it quasisymmetrically three point homogeneous} if for every pair of triples of distinct points $x,y,z \in X$ and $x',y',z' \in X$ there exists a quasisymmetric homeomorphism $f:X \to X$ so that $f(x) = x'$, $f(y) = y'$, and $f(z) = z'$. 

Chord-arc surfaces, first defined by Semmes in \cite{Semmes-CA1} and \cite{Semmes-CA2} are, roughly speaking, smooth, codimension-$1$ surfaces that come with quantitative control on the oscillation of the unit normal vector. More precisely, suppose that $\Gamma$ is a smooth connected hypersurface in $\mathbb{R}^d$ and that $\Gamma\cup\{\infty\}$ is smoothly embedded in $\mathbb{R}^{d+1}\cup\{\infty\}=\mathbb{S}^{d+1}$. Write $n(x)$ for the unit normal vector of $\Gamma$ at $x$ and $n_{z,r}$ for the average of $n$ over the ball $B(z,r)$ (with respect to the surface measure on $\Gamma$). For $\Gamma$ to be a \emph{chord-arc surface}, we ask that there is a constant $\gamma>0$ that bounds the BMO norm $\|n\|_*$ of $n$, i.e.,
$$
\sup_{x\in\Gamma, r>0} \frac{1}{|B(x,r)\cap \Gamma|} \int_{B(x,r)\cap \Gamma} |n - n_{z,r}| \leq \gamma,
$$
and satisfies the following condition:
$$ 
|\langle z-y, n_{z,r} \rangle | \leq \gamma r \text{ whenever } z,y\in \Gamma \text{ and } y\in B(z,r).
$$
A \emph{chord-arc surface with small constant} (sometimes `CASSC') is one for which the constant $\gamma$ is sufficiently small.

\subsection{Commentary}\label{subsec:qc-uniformization-parametrization-commentary}

This initial set of problems focuses on the interplay between deformations of measures and metrics in Euclidean space, quasiconformal and quasisymmetric parameterization of topological manifolds, and weak regularity properties for Euclidean hypersurfaces (e.g. the chord-arc condition). The story of the interactions between these concepts is closely linked to the origin of the subject of analysis in metric spaces, and work done throughout the 1990s on such questions initiated connections of the subject to other areas such as geometric topology. Such connections will be the focus of later sections of this survey. For these reasons, we devote some time here to reviewing the state of knowledge on this collection of ideas at the time when the Heinonen--Semmes problem list was produced.

The Jacobian of any quasiconformal map $f:\R^n \to \R^n$, $n \ge 2$, is a strong $A_\infty$ weight. This can be seen, for instance, by appealing to the quasisymmetry condition \eqref{eq:qs} to deduce that for any $x,y \in \R^n$, $|f(x)-f(y)| \simeq L_f(x,|x-y|) \simeq \ell_f(x,|x-y|) \simeq \cL^n(f(B_{x,y}))^{1/n}$ where $B_{x,y}$ is defined as in \eqref{eq:B-x-y}, $L_f(x,r) = \max\{|f(x)-f(a)|:|x-a|=r\}$, and $\ell_f(x,r) = \min\{|f(x)-f(b)|:|x-b|=r\}$. Moreover, it follows from a celebrated theorem of Gehring \cite{geh:higher-integrability} that strong $A_\infty$ weights are $A_\infty$ weights in the sense of Muckenhoupt. 

The process of deforming the geometry of Euclidean space $(\R^n,d_E)$ via a metric doubling measure $d\nu = w \, dx$, yielding the deformed metric space $(\R^n,d_w')$ where $d_w'$ satisfies \eqref{eq:dwprime}, was introduced by David and Semmes \cite{ds:strong-a-infinity-sobolev-quasiconformal}. In this paper, the authors showed that the analytic properties of the resulting metric measure space $(\R^n,d_w',w\,d\cL^n)$ closely resembled those of the undeformed space $(\R^n,d_E,d\cL^n)$. In particular, Sobolev and Poincar\'e inequalities hold for such weighted and deformed spaces with exponents identical to those in the undeformed geometry. These facts naturally led to the question of whether such deformed spaces could be realized, up to a bi-Lipschitz mapping, as subsets of a larger-dimensional Euclidean space. 

In \cite{sem:BL-strong-a-infinity}, Semmes posed the question to characterize strong $A_\infty$ weights $w$ on $\R^n$ and to determine whether or not every such weight is comparable to the Jacobian of a quasiconformal map. He also observed, \cite[Section 5]{sem:BL-strong-a-infinity}, that every $A_1$ weight in $\R^n$ is a strong $A_\infty$ weight. In a related paper \cite{sem:nonexistence}, he provided examples of strong $A_\infty$ weights $w$ on $\R^3$ which were not comparable to the Jacobian of any quasiconformal map of $\R^3$, which in turn meant that the deformed metric space $(\R^3,d_w')$ fails to be bi-Lipschitz equivalent to $\R^3$. Moreover, he gave examples of strong $A_\infty$ weights $w$ on $\R^n$ for some $n$ with the property that the deformed space $(\R^n,d'_w)$ fails to admit any bi-Lipschitz embedding into any $\R^N$. This circle of ideas was also seen to be related to the parameterizability question for chord arc surfaces with small constant. In particular, in \cite{Semmes-CA2} Semmes showed that chord arc surfaces in $\R^3$ with small constant admit quasisymmetric parameterizations, and in \cite{sem:BL-strong-a-infinity} he observed that control on the BMO norm of the normal vector $n_{z,r}$ could be leveraged to obtain good control on an associated strong $A_\infty$ weight.

Consequently, it was already known by the time of the preparation of \cite{hs:questions} that a positive answer to Question \ref{q:4} would imply a corresponding positive answer to Question \ref{q:1} and in turn to Question \ref{q:2}. In turn, a positive answer to both Questions \ref{q:1} and \ref{q:3} would also resolve Question \ref{q:4} in the affirmative. It was also known that Question \ref{q:3} has a positive answer in case the metric $d$ is a smooth Riemannian metric on $\Sph^2$, see \cite{hk:definitions} and \cite{ds:quantitative-rectifiability-lipschitz-mappings} for details. The restriction to the planar case is necessary. In \cite{sem:good} Semmes gave examples in dimensions $n\ge 3$ showing that the answer to the analog of Question \ref{q:3} could be `no', and in \cite{sem:selecta} (for $n \ge 5$) and \cite{sem:nonexistence} (for $n \ge 3$) he gave examples showing that the answer to the analog of Question \ref{q:4} could be `no'.

Semmes' example in $\R^3$, see \cite{sem:nonexistence}, involved a geometric realization of a classical object of geometric topology, {\it Antoine's necklace}. Specifically, he constructs a self-similar and totally disconnected subset $E \subset \R^3$, of Cantor set type, whose complement is not simply connected. The set $E$ is constructed as a decreasing union of compact sets $E_n$, where $E_n$ is a finite disjoint union of solid tori enjoying a suitable nontrivial linking property. It is well known that the non-simple connectivity of $\R^3 \setminus E$ implies a quantitative lower bound on the Hausdorff dimension $\dim_H E$. The required strong $A_\infty$ weight is chosen of the form $w(x) = \min \{1,\dist(x,E)^s \}$ for a suitable choice of $s>0$. The fact that the deformed space $(\R^3,d_w')$ is not bi-Lipschitz equivalent to $\R^3$ arises as a consequence of the aforementioned lower bound on $\dim_H E$ and a direct computation of an upper bound for the Hausdorff dimension of $E$ in the deformed metric $d_w'$. 

The topological wildness of Antoine's necklace features prominently in the preceding argument, and for some time it was not known whether analogous examples could be constructed in the plane. Questions \ref{q:1} and \ref{q:4} address this point. In 2002, Laakso introduced a family of planar metric graphs and an associated family of metric spaces which now go by the name {\it Laakso spaces}, which enabled him to resolve both Question \ref{q:1} and Question \ref{q:4} in the negative. Laakso's spaces are constructed using certain series-parallel and self-similar metric graphs. Equipped with a suitable path metric, such graphs become Ahlfors regular metric spaces supporting a Poincar\'e inequality. Alternatively, such a graph can be embedded as a suitable fractal subset $E$ of the usual plane $\R^2$. In \cite{laa:a-infinity}, Laakso considers $\R^2$ equipped with the deformed metric $d_w'$ defined as before, via a weight $w(x)$ given by a power of $\dist(x,E)$. He then proves the following theorem, which shows that even a substantially weaker version of Question \ref{q:4} also has a negative answer.

\begin{theorem}[Laakso]\label{th:laakso-1}
Let $G$ be a metric graph of the above type, and let $E \subset \R^2$ be the corresponding fractal realization. For suitable powers $s>0$, the weight $w(x) = \dist(x,E)^s$ defines a strong $A_\infty$ weight whose deformed metric geometry $(\R^2,d_w')$ does not admit a bi-Lipschitz embedding into any uniformly convex Banach space.
\end{theorem}

A fortiori, in view of the previously discussed results of David and Semmes, Theorem \ref{th:laakso-1} implies the following consequence.

\begin{theorem}[Laakso]\label{th:laakso-2}
None of the strong $A_\infty$ weights considered in Theorem \ref{th:laakso-1} are comparable to the Jacobian of a quasiconformal self-map of $\R^2$.
\end{theorem}

Although graphs of the type which feature in Theorem \ref{th:laakso-1} are nowadays commonly referred to as {\it Laakso graphs}, we note that the actual graph appearing in \cite{laa:a-infinity} is slightly different from the one which is now known as the Laakso graph; the latter was first introduced (to the best of our knowledge) in Lang and Plaut's influential paper \cite[Theorem 2.3]{LangPlaut}. The authors of \cite{LangPlaut} acknowledge that their construction is a variation on the examples considered in \cite{laa:a-infinity} and \cite{laa:PI}.

The Laakso graph is nowadays a standard example in analysis in metric spaces and metric embedding theory, arising in a host of different settings. A very brief, and woefully incomplete, list of references includes \cite{lmn:snowflakes}, \cite{tys:hyperspaces}, \cite{ck:inverse-limits}, \cite{ck:graphs}, \cite{li:markov}, \cite{ds:TST}, \cite{dko:laakso}, \cite{gar:markov}, and \cite{gar:thick}.

Laakso's construction left open Question \ref{q:2}, which in turn was answered by Bishop in \cite{bis:a-1}. Bishop's example consists of a carefully chosen non-self-similar Sierpi\'nski carpet $E \subset \R^2$, together with a weight $w$ on $\R^2$ such that $w$ is bounded away from zero on all of $\R^2$ and
\begin{equation}\label{eq:w-to-infinity}
\lim_{x \to E} w(x) = +\infty.
\end{equation}
The key feature of his example is that the image of $E$ under any quasiconformal map $f:\R^2 \to \R^2$ whose Jacobian is comparable to $w$, must necessarily contain a rectifiable curve $\gamma$. However, \eqref{eq:w-to-infinity} implies that $\det Dg = 0$ on $f(E)$, where $g = f^{-1}$, which in turn would imply that $g$ is constant on $\gamma$. The latter conclusion is incompatible with the assumption that $f$ (and hence also $g$) is a homeomorphism. It follows that $w$ cannot be comparable to the Jacobian of any planar quasiconformal map. Bishop also observes that his example, in conjunction with the previous observations of Semmes regarding geometric deformations via $A_1$ weights, generates a deformed space which bi-Lipschitz embeds into some finite-dimensional Euclidean space. By a direct construction using the explicit example under consideration, he in fact showed that his example embeds bi-Lipschitzly into $\R^3$. In other words, he establishes the following result.

\begin{theorem}[Bishop]
There exists a linearly locally connected and Ahlfors $2$-regular surface $S \subset \R^3$ which is not bi-Lipschitz equivalent to $\R^2$.
\end{theorem}

Question \ref{q:3} asks for a weaker notion of parametrization, under the same assumptions as in Question \ref{q:4}. Questions \ref{q:3} and \ref{q:4} are instances of a general uniformization program in metric analysis, namely, the search for intrinsic properties sufficient for equivalence with a suitable model space. Such questions are of interest in both the bi-Lipschitz and the quasisymmetric category. Bi-Lipschitz and quasisymmetric characterizations of the unit circle $\Sph^1$ were already known by the time that \cite{hs:questions} was published. A metric space $(X,d)$, assumed {\it a priori} to be homeomorphic to $\Sph^1$, is quasisymmetrically equivalent to $\Sph^1$ if and only if it is doubling and satisfies the {\it bounded turning property}. Moreover, such $(X,d)$ is bi-Lipschitz equivalent to $\Sph^1$ if and only if it is doubling and satisfies the {\it chord-arc property}.\footnote{A topological circle $(X,d)$ is {\it bounded turning} if there exists $C>0$ so that $\min\{\diam(\gamma_1),\diam(\gamma_2)\} \le C d(x,y)$ for all distinct points $x,y \in X$, where $\gamma_1$ and $\gamma_2$ denote the two connected components of $X \setminus \{x,y\}$. Similarly, $(X,d)$ is {\it chord-arc} if there exists $C>0$ so that $\min\{\length(\gamma_1),\length(\gamma_2)\} \le C d(x,y)$ for all $x\ne y$. Note that the chord-arc assumption implies that any pair of points in $X$ can be joined by a rectifiable curve.} We refer the reader to \cite[Section 3]{bonk:icm} for a broad overview of the quasisymmetric uniformization program and for more on these results. Questions \ref{q:3} and \ref{q:4} concern similar characterizations for the $2$-sphere $\Sph^2$.

While not explicitly discussed in \cite{hs:questions}, substantial motivation for the problem of quasisymmetric uniformization of the $2$-sphere $\Sph^2$ comes from geometric group theory and Cannon's conjecture. A {\it Gromov hyperbolic group} is a finitely generated group $G$ whose Cayley graph $C(G)$ (equipped with the geodesic path metric such that every edge has length one) is $\delta$-hyperbolic in the sense of Gromov, for some $\delta>0$. Recall that a geodesic metric space $(X,d)$ is said to be {\it $\delta$-hyperbolic} if geodesic triangles satisfy the `thin triangles' condition: for any points $x,y,z \in X$ and any choice of geodesics $[x,y]$, $[y,z]$, and $[x,z]$, it holds true that any of the edges $[x,y]$, $[y,z]$, and $[x,z]$ is contained in the metric $\delta$-neighborhood of the other two edges. To each Gromov hyperbolic geodesic metric space $(X,d)$, and in particular to each Gromov hyperbolic group (viewed as its Cayley graph with geodesic metric), one associates a notion of boundary at infinity. This space, the {\it Gromov boundary} $\partial_\infty X$, consists of all infinite geodesic rays emanating from a fixed basepoint $x_0 \in X$ under a suitable equivalence relation. 

The {\it Gromov product} of two points $x,y$ in a Gromov hyperbolic group $G$ with respect to a basepoint $x_0 \in G$ is defined to be $(x|y)_{x_0} = \tfrac12(d(x,x_0)+d(y,x_0)-d(x,y))$. Taking into account the defining equivalence relation, one observes that the Gromov product extends to pairs of points $\xi,\eta$ in $\partial_\infty G$. Following a procedure closely related to the David--Semmes deformations discussed above, one defines a family of metrics $d_{x_0,\eps}$ on $\partial_\infty G$. More precisely, there exists $\eps_0 = \eps_0(\delta) > 0$ depending on the hyperbolicity parameter $\delta$ so that for every $x_0 \in G$ and every $0<\eps<\eps_0$, there exists a metric $d_{x_0,\eps}$ on $\partial_\infty G$ and a constant $L \ge 1$ so that
$$
\frac1L d_{x_0,\eps}(\xi,\eta) \le e^{-\eps (\xi|\eta)_{x_0}} \le L d_{x_0,\eps}(\xi,\eta), \qquad \forall\,\xi,\eta \in \partial_\infty G.
$$
Such metrics on $\partial_\infty G$ are called {\it visual metrics}. It is not difficult to check that any two visual metrics on $\partial_\infty G$ are quasisymmetrically equivalent, i.e., the identity map between these two metrics is a quasisymmetric homeomorphism. Consequently, quasisymmetric invariants of the Gromov boundary equipped with any such visual metric are of interest for understanding the large-scale geometry of $G$ and (ultimately) algebraic properties of $G$. For instance, if $\partial_\infty G$ is homeomorphic to the Cantor set $C$, then $\partial_\infty G$ is quasisymmetrically equivalent to $C$. In this case, the group $G$ is {\it virtually free}, i.e., it is a free group up to a finite index subgroup. Similarly, if $\partial_\infty G$ is homeomorphic to $\Sph^1$, then $\partial_\infty G$ is quasisymmetrically equivalent to $\Sph^1$ and $G$ is {\it virtually Fuchsian}, i.e., up to a finite index subgroup it can be realized as a group of discrete, cocompact isometries of the real hyperbolic plane $H^2_\R$. The two-dimensional analog is a celebrated conjecture of Jim Cannon.

\begin{conjecture}[Cannon]\label{conj:cannon}
Let $G$ be a Gromov hyperbolic group such that $\partial_\infty G$ is homeomorphic to $\Sph^2$. Then up to a finite index subgroup, $G$ can be realized as a group of discrete, cocompact isometries of $H^3_\R$. That is, $G$ is {\it virtually Kleinian}.
\end{conjecture}

A standard argument shows that the conclusion in Conjecture \ref{conj:cannon} is equivalent to showing that $\partial_\infty G$ is quasisymmetrically equivalent to the round sphere $\Sph^2$. One is naturally led to inquire what are appropriate metric conditions on a topological $2$-sphere which suffice to conclude quasisymmetric equivalence.

We also note in passing that the analogous three-dimensional statement is false. Any discrete and cocompact group of isometries of the complex hyperbolic plane $H^2_\C$ yields an example of a Gromov hyperbolic group $G$ whose boundary $\partial_\infty G$ is homeomorphic to $\Sph^3$, but is not quasisymmetrically equivalent to $\Sph^3$. In fact, such a boundary is quasisymmetrically equivalent to a suitable one-point conformal compactification of the sub-Riemannian Heisenberg group. The Heisenberg group will be discussed further below in section \ref{sec:heisenberg}.

In a series of papers \cite{bk:spheres, bk:rigidity, bk:rigidity2, bk:conformal} Bonk and Kleiner studied the quasisymmetric uniformization problem for $\Sph^2$ and obtained strong partial results in support of Cannon's conjecture. In \cite{bk:spheres}, Bonk and Kleiner answered Question \ref{q:3} in the affirmative.

\begin{theorem}[Bonk--Kleiner]\label{thm:bk}
Let $(X,d)$ be a metric space which is homeomorphic to $\Sph^2$ and is both linearly locally contractible and Ahlfors $2$-regular. Then $(X,d)$ is quasisymmetrically equivalent to $\Sph^2$.
\end{theorem}

Remarkably, there are now four essentially distinct approaches to Theorem \ref{thm:bk}, each relying on a different set of techniques.
\begin{enumerate}[(1)]
\item The original proof by Bonk and Kleiner \cite{bk:spheres} relies on \textbf{circle packing}. In essence, Bonk and Kleiner triangulate the metric $2$-sphere $X$ in a controlled way and apply the classical circle packing theorem of Andreev--Koebe--Thurston to find a circle packing of the sphere with a combinatorially equivalent dual graph. This yields a coarse mapping to $\mathbb{S}^2$ whose quasisymmetric distortion can, with careful estimates, be controlled. The final quasisymmetry arrives by passing to a limit as the scale of the approximation tends to zero.
\item Kai Rajala \cite{raj:spheres} gave a new proof for Theorem \ref{thm:bk} which relies on properties of \textbf{harmonic functions on non-smooth surfaces}. Rajala in fact proves a more general result than Theorem \ref{thm:bk}, which provides a novel necessary and sufficient condition for a planar surface of locally finite two-dimensional Hausdorff measure to admit a \emph{geometrically quasiconformal} homeomorphism to a planar domain.
\item An argument by Lytchak and Wenger \cite{lw:spheres} proceeds via a study of \textbf{area-minimizing disks in metric spaces}. These authors show that if $X$ is an Ahlfors $2$-regular, linear locally contractible metric space homeomorphic to the closed unit disk $\mathbb{D}$ in the plane, then $X$ admits a parametrization of minimal `energy', and this parametrization is quasisymmetric. The notion of `energy' here is a version of the $W^{1,2}$ Sobolev energy for metric space valued maps called the `Reshetnyak energy'. (Roughly speaking, it is the $L^2$-norm of the minimal weak upper gradient of the map; see subsection \ref{subsec:mms-definitions-and-notation}.) Later, Meier and Wenger \cite{meierwenger} combined energy minimization methods with some more flexible techniques to prove a more general parametrization result that also recovers Theorem \ref{thm:bk} in yet a different way.
\item Finally, Ntalampekos and Romney \cite{nr:poly} give a different proof which makes use of \textbf{polyhedral approximation results for non-smooth surfaces}. They show that a metric surface $X$ of finite $2$-dimensional Hausdorff measure can be approximated by polyhedral surfaces in a useful way. From these approximations, they construct a surjective map from a smooth surface to $X$ with a certain weak, one-sided, quasiconformality property. Under the assumptions of Theorem \ref{thm:bk}, this property can be upgraded to true quasisymmetry.
\end{enumerate}
A more detailed overview of these approaches can be found in the recent survey \cite{nta:survey}. 

As a consequence of Theorem \ref{thm:bk}, and by making use of results of Keith--Laakso \cite{kl:conformal-assouad-dimension} on the structure of metric spaces minimal for conformal Assouad dimension, Bonk and Kleiner showed in \cite{bk:conformal} that the conclusion of Conjecture \ref{conj:cannon} holds for a hyperbolic group $G$ with $2$-sphere boundary if and only if the so-called \emph{Ahlfors regular conformal dimension} of $\partial_\infty G$, namely, the quantity
$$ 
\inf\{ Q : \partial_\infty G \text{ is quasisymmetric to an Ahlfors } Q\text{-regular space} \}
$$
is attained as a minimum. The attainment question for conformal dimension remains quite subtle, and we refer the reader to Section 6 of \cite{bk:conformal} and \cite{mt:cdim} for further discussion of conformal dimension, which has developed into a major area of research.

\smallskip

Questions \ref{q:5} and \ref{q:6} remain open in general. Under a stronger set of hypotheses (namely, that the indicated homogeneity is obtained via a family of uniformly quasi-M\"obius homeomorphisms), an affirmative answer to Question \ref{q:5} was provided by Bonk and Kleiner in \cite{bk:rigidity}.

\section{Bi-Lipschitz parametrizations and embeddings}\label{sec:bl-parametrizations-embeddings}
\begin{question}\label{q:7}
\it Do smooth $n$-dimensional chord arc surfaces with small constant in $\R^{n+1}$, $n \ge 2$, admit bi-Lipschitz parametrizations?
\end{question}

This question remains open. 

\begin{question}\label{q:8}
\it If an Ahlfors regular metric space admits a regular map into some Euclidean space $\mathbb{R}^d$, then does it admit a bi-Lipschitz map into another, possibly different, Euclidean space?
\end{question}

The answer is NO if $d\geq 2$, by work of David and Eriksson-Bique \cite{deb:slit-carpets}, and remains open for $d=1$.

\begin{question}\label{q:9}
\it If an Ahlfors $n$-regular metric space has big pieces of Lipschitz images of $\R^n$, $n \ge 2$, is it then uniformly rectifiable of dimension $n$?
\end{question}

The answer is YES, as proved by Schul \cite{sch:decomposition}. 

\subsection{Definitions and notation}\label{subsec:bl-parametrizations-embeddings-definitions-and-notation}

Intermediate between general Lipschitz mappings and bi-Lipschitz embeddings lie the (David--Semmes) regular mappings, first introduced in \cite{David_operateurs} in a different form and given their current definition in \cite{ds:fractured-fractals}, \cite{ds:regular}. A map $g\colon (X,d) \rightarrow (Y,d')$ is \emph{regular} if it is Lipschitz and there is a constant $C>0$ such that for each ball $B(y,r)\subseteq Y$, the set $f^{-1}(B(y,r))\subseteq X$ can be covered by at most $C$ balls of radius $Cr$. Regular mappings need not be injective, but they permit only a controlled amount of folding; in particular, a regular mapping with constant $C$ is at most $C$-to-$1$.

One of the most important classes of objects in geometric measure theory are $n$-\emph{rectifiable sets}, or more generally $n$-\emph{rectifiable metric spaces}. These are spaces that are covered, up to zero $\mathcal{H}^n$-measure, by Lipschitz images of subsets of Euclidean space; see \cite{Mattila}, \cite{Federer} for more precise definitions.

In \cite{ds:uniform-rectifiability}, David and Semmes pioneered the study of a quantitative form of rectifiability called uniform rectifiability. An Ahlfors $n$-regular metric space $X$ is \emph{uniformly $n$-rectifiable}, with constants $\theta, M$, if each ball $B(x,r)$ in $X$, $0<r\leq \diam(X)$, contains a subset $E$ with $\cH^n(E)\geq \theta r^n$ that is $M$-bi-Lipschitz to a subset of $\mathbb{R}^n$. This condition is much stronger than classical rectifiability and has close connections to the study of singular integrals.

An \textit{a priori} weaker quantitative condition on $X$ is that it contains \emph{big pieces of Lipschitz images of $\R^n$}, with constants $\rho, L$. This means that each ball $B(x,r)$ in $X$, $0<r\leq \diam(X)$, contains a subset $F$ with $\mathcal{H}^n(F)\geq \rho r^n$ that is an $L$-Lipschitz image of a subset of $B(0,r)\subseteq\mathbb{R}^n$.\footnote{Note that the definition given on \cite[p.\ 4]{hs:questions} is not quite correct.}

\subsection{Commentary}\label{subsec:bl-parametrizations-embeddings-commentary}
Questions \ref{q:6} and \ref{q:7} originate in work of Semmes, which in turn stem from the search for refinements to the classical `topological disk' theorem of Reifenberg \cite{Reifenberg}. We state Reifenberg's theorem as presented in \cite{dt:reifenberg2}:
\begin{theorem}[Reifenberg, 1960]\label{thm:reifenberg}
For each $1\leq d \leq n$ and $\tau\in (0,\frac{1}{10})$, there is a parameter $\epsilon>0$ with the following property: Let $E$ be a closed set in $\mathbb{R}^n$ containing the origin. Suppose that for each $x\in E$ and $0\leq r \leq 10$, there is an affine $d$-plane $P_{x,r}$ through $x$ such that
$$ \dist(y,P_{x,r}) \leq \epsilon r \text{ and } \dist(z, E) \leq \epsilon r$$
whenever $y\in E\cap B(x,r)$ or $z\in P_{x,r} \cap B(x,r)$.

Then there is a bijection $g\colon \mathbb{R}^n \rightarrow \mathbb{R}^n$ such that
$$ E\cap B(0,1) = g(\mathbb{R}^d \times \{0_{\mathbb{R}^{n-d}}\}) \cap B(0,1)$$
and
$$ \frac14 |x-y|^{1+\tau} \leq |g(x)-g(y)| \leq 3|x-y|^{1-\tau} \text{ for all } x,y\in\mathbb{R}^n \text{ with } |x-y|\leq 1.$$
\end{theorem}
Reifenberg's theorem tells us that a set that looks `flat' with a (small) uniform error at every location and scale must be a topological disk; in fact, it must be a bi-H\"older disk with exponent close to $1$. Much work in quantitative rectifiability has been devoted to refining or extending Reifenberg's theorem in various directions, for instance, to allow for different notions of `flatness' or stronger types of parametrization. The work of David and Toro in \cite{dt:reifenberg, dt:reifenberg2} provides major examples.

Questions \ref{q:6} and \ref{q:7} are in some sense part of this program. Of course, a positive answer to Question \ref{q:7} would imply a positive answer to Question \ref{q:6}, though Semmes conjectures in \cite{Semmes-CA2} that the answer to both of these questions is `no'. Semmes himself provides some partial results in \cite{Semmes-CA2}. Later work of Toro \cite{tor:bilipschitz}, David--Toro \cite{dt:reifenberg2}, and Merhej \cite{mer:geometry,mer:parameterizations} provides bi-Lipschitz parametrizations under stronger quantitative conditions on the normal vector.

\begin{figure}[h]
\includegraphics[width=8cm]{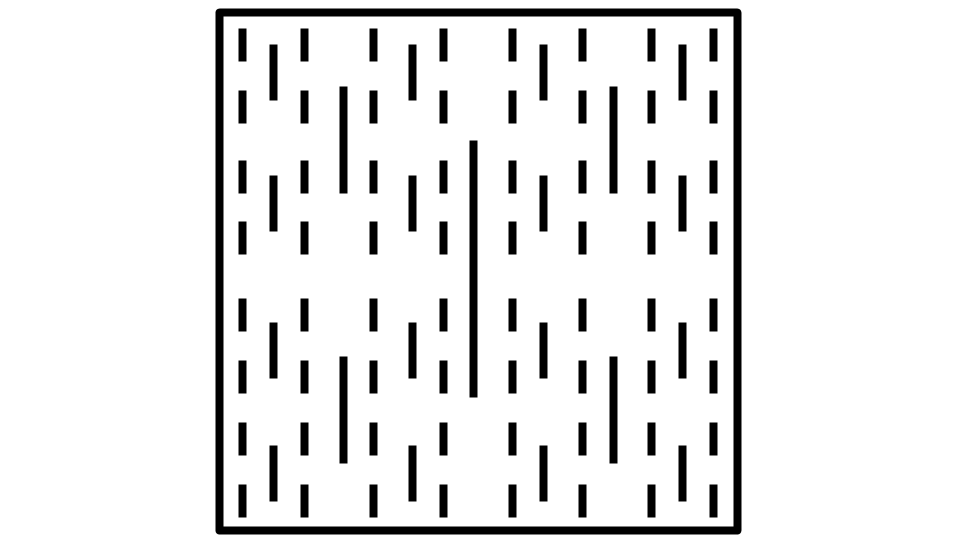}
\centering
\caption{An approximation to the slit carpet $M$.}
\label{fig:slit}
\end{figure}

The counterexample that yields a negative answer to Question \ref{q:8} in \cite{deb:slit-carpets} for $d=2$ (from which the $d\geq 2$ case follows immediately) is the so-called {\it slit carpet}, studied extensively by Merenkov in a different context \cite{mer:co-hopfian}. To build the slit carpet, informally, start with the unit square $[0,1]^2$ and remove a vertical ``slit'' of length $\frac12$ through the center. Then repeat this, scaled down, in each dyadic sub-cube of $[0,1]^2$; see Figure \ref{fig:slit}. Equip the complement of these slits with its \emph{intrinsic path metric}, and take the completion. The resulting space $M$ is the slit carpet, so-called because it is homeomorphic to the classical Sierpi\'nski carpet. (For a more careful definition, see \cite{mer:co-hopfian} or \cite{deb:slit-carpets}.)

Merenkov observed that $M$ is a compact, Ahlfors $2$-regular space equipped with a David--Semmes regular map into $\R^2$ given by `collapsing the slits'. In \cite{deb:slit-carpets}, it is shown that $M$ admits no bi-Lipschitz embedding into any $\R^N$, or even any uniformly convex Banach space. An alternative proof of the non-embedding into $\R^N$ is given in \cite{deb:splitting}.

Question \ref{q:8} is a special case of the `Bi-Lipschitz Embedding Problem': {\it Find geometric conditions for a metric space to admit (or fail to admit) a bi-Lipschitz embedding into some finite-dimensional Euclidean space.} There are no simple necessary and sufficient conditions, and the problem has driven much work in analysis, geometry, and theoretical computer science; see for instance \cite{naor:icm}. Digressing slightly, we mention that perhaps the most fundamental open problem in the area is the `Lang--Plaut Problem' \cite{LangPlaut}: {\it Does every doubling subset of a Hilbert space admit a bi-Lipschitz embedding into a finite-dimensional Euclidean space?} While the Lang--Plaut problem remains open in full generality, an interesting partial result by Lafforgue and Naor \cite{ln:lp} states that for each $p>2$, $L^p$ contains a doubling subset that does not admit any bi-Lipschitz embedding into any finite-dimensional Euclidean space.

The David--Semmes theory of uniform rectifiability \cite{ds:uniform-rectifiability} is extremely well-developed for subsets of Euclidean space. Historically, much less has been known in metric spaces, though there has been significant recent progress; see, e.g., \cite{naoryoung1, naoryoung2, batehydeschul}. Question \ref{q:9} was an early strike in this direction. Like Questions \ref{q:22}--\ref{q:24} below, it asks whether a known property of Lipschitz mappings between Euclidean spaces holds in more general metric spaces. Here the property under consideration is a type of `quantitative inverse function theorem' for Lipschitz mappings: a Lipschitz mapping from $\mathbb{R}^n$ into a metric space admits a Lipschitz inverse, with quantitative bounds, on a definite fraction of its image. Such a result was already known to hold in a qualitative sense by a result of Kirchheim \cite{Kirchheim}, or quantitatively under the additional strong assumption that the image was a subset of Euclidean space (see \cite{David88}, \cite{jones:bilipschitz}, \cite{ds:uniform-rectifiability}).
The innovation of Schul in \cite{sch:decomposition} was to provide a purely metric approach to this type of quantitative theorem. We give a precise statement of a slightly special case of Schul's theorem.

\begin{theorem}[Schul]\label{thm:schul}
Fix $\alpha\in (0,1)$, and $n\geq 1$. There are constants $M=M(\alpha,n)$ and $c=c(n)$ with the following property: If $f\colon [0,1]^n \rightarrow X$ is a $1$-Lipschitz mapping into a metric space, then there are sets $F_1, \dots, F_M \subseteq [0,1]^n$ such that
\begin{enumerate}[(i)]
\item if $1\leq i\leq M$ and $x,y\in F_i$, then
\begin{equation}\label{eq:schul}
\alpha|x-y| \leq d(f(x),f(y)) \leq |x-y|,
\end{equation}
and
\item $\mathcal{H}^n_\infty\left(f([0,1]^n \setminus \cup_{i=1}^M F_i)\right) \leq c\alpha$.
\end{enumerate}
\end{theorem}

Here $\mathcal{H}^n_\infty$ denotes the \emph{$n$-dimensional Hausdorff content}:
$$ \mathcal{H}^n_\infty(A) = \inf\{ \sum \diam(U_i)^n : \{ U_i\} \text{ is a cover of } A\}.$$
In this definition, unlike for the Hausdorff measure $\mathcal{H}^n$, there is no requirement that the sets $U_i$ in the cover are small. We note that using content rather than measure in (ii) is necessary for Theorem \ref{thm:schul}. Indeed, suppose $n=1$ is given, $c=c(n)$ is as in Theorem \ref{thm:schul}, $\alpha\ll c$, and $M=M(\alpha,n)$ is as in Theorem \ref{thm:schul}. Let $f\colon [0,1]\rightarrow X$ be the arc-length parametrization of a Jordan arc of length $1$ but diameter $\ll_{M,c} \alpha$.\footnote{Imagine a `ball of yarn''.} Then it is not hard to see that one cannot have sets $F_i$ satisfying (i) and (ii) in the theorem, if $\mathcal{H}^n_\infty$ is replaced by $\mathcal{H}^n$. In this case, (the correct version of) Theorem \ref{thm:schul} would simply set each $F_i=\emptyset$, which is permitted since the \emph{content} of the image of $f$ is small.

Theorem \ref{thm:schul} implies a positive answer to Question \ref{q:9} without much further work.\footnote{Schul (personal correspondence) remarks on the Heinonen--Semmes list, ``I did not look at the list of open problems and try to solve them.  Rather, it was a way of seeing that the world is full of unknowns that are easy to state.''} Further quantitative results on Lipschitz mappings into metric spaces, including `quantitative implicit function theorems' for Lipschitz mappings in a style inspired by Theorem \ref{thm:schul}, appear in \cite{as:metric-differentiation}, \cite{as:hard-sard}, \cite{ds:lipschitz-decompositions}. The main engine behind many of these arguments are `quantitative differentiation' results for Lipschitz mappings on $\mathbb{R}^n$, which say that Lipschitz mappings are `nearly affine' on most locations and scales, even coarse scales where classical differentiability does not apply. See  \cite{Dorronsoro}, \cite{as:metric-differentiation},
and \cite[Appendix B]{Gromov} for more.

Question \ref{q:9} and the discussion of the previous paragraph concern Lipschitz mappings from $\mathbb{R}^n$ into metric spaces. It is natural to ask if there is a similar `quantitative inverse function theorem' for Lipschitz mappings from other domains, for instance the Heisenberg group. Here the answer is known to be `no'; see Question \ref{q:24} and its solution below.

\section{Geometric topology}\label{sec:geometric-topology}

\begin{question}\label{q:10}
\it Does the space $W \times \R^k$ admit for some $k \ge 1$ a quasisymmetric map onto $\R^{3+k}$, where $W$ is a linearly locally contractible and Ahlfors $3$-regular geometric realization of the decomposition space $\R^3/\cW$ associated with the Whitehead continuum $\Wh$?
\end{question}

The answer is NO. By a result of Heinonen and Wu \cite{hw:whitehead}, there exists linearly locally contractible and Ahlfors $3$-regular geometric realizations $W$ of $\R^3/\cW$ for which the product spaces $W\times \R^k$ are not quasisymmetrically equivalent to the Euclidean space $\R^{3+k}$.

\begin{question}\label{q:11}
\it Does the space $B \times \R^k$ admit for some $k \ge 1$ a quasisymmetric map onto $\R^{3+k}$, where $B$ is a linearly locally contractible and Ahlfors $3$-regular geometric realization of the decomposition space $\R^3/\cB$ associated with the Bing double $\Bd$?
\end{question}

The answer is NO. By a work of Pankka and Wu \cite{pw:semmes}, there exist geometric realizations $B$ of $\R^3/\cB$, which are linearly locally contractbile and Ahlfors $3$-regular but for which $B\times \R^k$ is not quasisymmetrically equivalent to $\R^{3+k}$.

\begin{question}\label{q:12}
\it Is there a quasisymmetric map from a polyhedral Edwards sphere $X$ onto $\Sph^5$?
\end{question}

This question remains open.

\begin{question}\label{q:13}
\it Is there a homeomorphism from a polyhedral Edwards sphere $X$ onto $\Sph^5$ that preserves sets of Hausdorff $5$-measure zero?
\end{question}

This question remains open.

\begin{question}\label{q:14}
\it Is there a homeomorphism from a polyhedral Edwards sphere $X$ onto $\Sph^5$ that belongs to the Sobolev space $W^{1,p}(X)$ for some $p \ge 1$?
\end{question}

This question remains open.

\subsection{Definitions and notation}\label{subsec:geometric-topology-definitions-and-notation}

A {\it decomposition space} for Euclidean space $\R^3$ is 
a quotient space $\R^3/\mathcal{P}$, that is, a partition $\mathcal P$ of $\R^3$ with a quotient topology; the partition $\mathcal P$ is called a \emph{decomposition of $\R^3$}.
For the purposes of this survey, we define a {\it geometric realization} of a decomposition space $\R^3/\mathcal{P}$ to be a closed set $F \subset \R^4$ so that $F$ (equipped with the restriction of the Euclidean metric of $\R^4$) is homeomorphic to $\R^3/\mathcal{P}$.

The {\it Whitehead continuum} is a closed subset $\Wh \subset \R^3$ obtained by the following iterative procedure. Begin with a closed solid torus $T_0 \subset \R^3$, i.e.~a compact subset homeomorphic to $\bar B^2 \times \mathbb S^1$, and consider a second solid torus $T_1$ contained in the interior of $T_0$ so that $T_1$ is homotopically trivial but isotopically nontrivial; see Whitehead \cite{Whitehead} for the precise construction or e.g.~Daverman \cite[Figure 9-7]{DavermanR:Decm} for an illustration.
Next, fix a homeomorphism $h$ from $T_0$ onto $T_1$ and define $T_2 := h(T_1)$. Inductively, set $T_{m+1} := h(T_m)$ for each $m \ge 1$. Finally, let $\Wh = \cap_{m=1}^\infty T_m$. The associated decomposition space $\R^3/\cW$ consists of the partition $\mathcal{W}$ whose elements are the singleton sets $\{x\}$ for each $x \not\in \Wh$, and the set $\Wh$ itself.

The {\it Bing double} $\Bd \subset \R^3$ is defined in a similar fashion, again starting with a closed solid torus $\Sigma_0 = T_0$. The second body $\Sigma_1$ is assumed to be contained in the interior of $\Sigma_0$ and to consist of the union of two disjoint closed solid tori $\Sigma_1 := T_1^1 \cup T_1^2$, so that each $T_1^i$, $i=1,2$, is isotopically trivial in $\Sigma_0$, while the union $T_1^1 \cup T_1^2$ is homotopically trivial but isotopically nontrivial in $\Sigma_0$; see Bing \cite{bing:double} for the precise construction or e.g.~Daverman \cite[Figure 9-1]{DavermanR:Decm} for an illustration. Next, for $i=1,2$ fix homeomorphisms $h^i$ from $T_0$ to $T_1^i$ and define $T_2^{ij} := h^i(T_1^j)$ and $\Sigma_2 = \cup_{i,j} T_2^{ij}$. Inductively, define $T_m^w$ for each $m \ge 1$ and each word $w \in \{1,2\}^m$ and set $\Sigma_m = \cup_w T_m^w$. Finally, let $\Bd = \cap_{m=1}^\infty \Sigma_m$. The associated decomposition space $\R^3/\cB$ consists of the partition $\mathcal{B}$, whose elements are singleton sets $\{x\}$ for each $x \not\in \Bd$, and components of the set $\Bd$.

The {\it Edwards sphere} (or {\it Cannon--Edwards sphere}) is the double suspension $X = \Sigma^2 H^3$ of a space $H^3$, which is a homology $3$-sphere but not a topological $3$-sphere. An instance of such a homology $3$-sphere is the Poincar\'e homology sphere (Poincar\'e dodecahedral space). Polyhedral models for such homology $3$-spheres can be constructed, in which case the associated double suspension $\Sigma^2 H^3$ can also be given a polyhedral structure. The Cannon--Edwards double suspension theorem asserts that in this setting, $\Sigma^2 H^3$ is a topological $5$-sphere. 

For a polyhedral Edwards $5$-sphere $X = \Sigma^2 H^3$, the Sobolev space $W^{1,p}(X)$ is defined to consist of mappings $f:X \to \Sph^5 \subset \R^6$ such that $f|_\Delta \in W^{1,p}(\Delta)$ for each $5$-simplex $\Delta$ contained in $X$. Note that any such simplex $\Delta$ may be identified with the closure of a domain in $\R^5$, and the Sobolev space consists of continuous mappings $f|_\Delta : \Delta \to \R^6$ so that $f$ lies in $L^p$ and $f$ has weak partial derivatives which also lie in $L^p$.

\subsection{Commentary}\label{subsec:geometric-topology-commentary}

Questions \ref{q:10}--\ref{q:14} are geometric counterparts of classical topological questions in geometric topology on exotic factors of Euclidean spaces, wild homeomorphisms, and the manifold recognition problem. 

A topological space $X$ is an \emph{exotic factor} if $X$ is not homeomorphic to a Euclidean space $\R^n$ but there exists $k>0$ for which $X\times \R^k$ is homeomorphic to a Euclidean space. The space $X\times \R^k$ is sometimes called a \emph{stabilization of $X$ by $\R^k$}.

By a classical theorem, obtained independently by Edwards and Miller \cite{Edwards-Miller} and Pixley and Eaton \cite{Pixley-Eaton}, \emph{for a closed-$0$-dimensional cell-like decomposition $G$ of $\mathbb S^3$, the space $(\mathbb S^3/G) \times \R^1$ is homeomorphic to $\mathbb S^3\times \R$}. As a consequence, the aforementioned decomposition space $\R^3/\cW$ in Question \ref{q:10} is an exotic factor, although it is not a manifold itself. The decomposition space $\R^3/\cB$ is, by a famous theorem of Bing \cite{bing:double}, homeomorphic to $\R^3$. The homeomorphism $\R^3/\cB \to \R^3$ constructed by Bing is a classical example of a so-called \emph{wild homeomorphism}. We refer to the monograph of Daverman \cite{DavermanR:Decm} for the role of decomposition spaces in geometric topology.

Questions \ref{q:10} and \ref{q:11} stem from the observation that, although a decomposition space $\R^3/\mathcal P$ for example in the case of exotic factors is metrizable, the quotient map $\R^3 \to \R^3/\mathcal P$ does not induce a canonical metric on $\R^3/\mathcal P$. Semmes showed in \cite{sem:good} that Bing's double admits a geometric realization $\R^3/\cB \hookrightarrow \R^4$ which maps the tori $T^w_m$ to similar round tori $\widetilde T^w_m$ in $\R^4$ having diameters tending geometrically to zero as the length of the word $w$ tends to infinity. The realization $\R^3/\cB \hookrightarrow \R^4$ can be chosen so that the metric $d_{\cB}$ induced by the embedding is linearly locally contractible and Ahlfors $3$-regular. 
Semmes proved that, despite these properties, the metric space $(\R^3/\cB,d_{\cB})$ -- or more concretely, the image of the embedding $\R^3/\cB \hookrightarrow \R^4$ -- is not quasisymmetric to $\R^3$. Semmes' embedding argument can easily be modified to give an embedding $\R^3/\cW \hookrightarrow \R^4$ of the Whitehead space $\R^3/\cW$ into $\R^4$ for which the induced geometric realization metric $d_{\cW}$ on $\R^3/\cW$ is again linearly locally contractible and Ahlfors $3$-regular. We refer to an essay of Wu \cite{Wu-Semmes-spaces} for a more detailed discussion on these ideas and related constructions of metrics.

Using the classical terminology in topology, Questions \ref{q:10} and \ref{q:11} may now be interpreted to ask \emph{whether spaces $W=(\R^3/\cW,d_\cW)$ and $B=(\R^3/\cB,d_\cB)$ are exotic quasisymmetric factors of some Euclidean space $\R^{3+k}$ for $k>0$}. The seminal result of Heinonen and Wu \cite{hw:whitehead} for the geometric realization $W$ of the Whitehead space $\R^3/\cW$, and its generalization by Pankka and Wu \cite{pw:semmes} for the geometric realization $B$ of the Bing double $\R^3/\cB$, show that this is not the case.

\begin{theorem}[Heinonen--Wu, Pankka--Wu]
Let $X$ be either a geometric realization $W$ of the Whitehead space $\R^3/\cW$ or a geometric realization $B$ of the Bing double $\R^3/\cB$. Then the product spaces $X\times \R^k$, for $k\ge 0$, are linearly locally contractible and Ahlfors $(3+k)$-regular Loewner spaces that are not quasisymmetrically equivalent to $\R^{3+k}$.
\end{theorem}

The non-existence of a quasisymmetric homeomorphism in this result is based on a homological intersection counting argument of Freedman and Skora \cite{Freedman-Skora}, also used by Semmes in \cite{sem:good}. It should be noted that, whereas in \cite{Freedman-Skora} and \cite{sem:good} the intersection counting is used to obtain lower bounds for lengths of paths, in \cite{hw:whitehead} and \cite{pw:semmes} intersection counting is used to estimate modulus of $(k+1)$-dimensional surface families. To our knowledge, the method of Heinonen and Wu in \cite{hw:whitehead} is one of the first instances where the modulus of a surface family is used effectively to obtain results about quasiconformal mappings. The same intersection counting argument for moduli of surface families can also be applied to other decomposition spaces that arise as partitions associated to defining sequences, that is, nested sequences of $3$-manifolds with boundary. We refer to \cite{pw:semmes} for a discussion on geometric realizations of more general decomposition spaces and on the interplay between the geometric parameters in the geometric realization metrics and the intersection counting data of the defining sequences.

Although there has been no progress in Questions \ref{q:12}--\ref{q:14}, we recall some background on the topological origin of these problems. The manifold recognition problem asks, which qualitative properties of a topological space yield a local parametrization with a Euclidean space; see e.g.~Cannon \cite{CannonJ:Recpwi}. 
This problem was settled independently by Edwards \cite{EdwardsR:Topmcl} (see also \cite{EdwardsR:Sushs}) and Cannon \cite{CannonJ:Shrcld} by results on suspensions of homology spheres.

In Questions  \ref{q:12}--\ref{q:14} on geometry and analysis of the Cannon--Edwards sphere, the starting point is the aforementioned observation that double suspension of a $3$-dimensional polyhedral structure of a homology sphere $H^3$ yields a $5$-dimensional polyhedral structure on $\Sigma^2 H^3$, which carries a natural polyhedral metric making $\Sigma^2 H^3$ a metric $5$-sphere. It was observed by Siebenmann and Sullivan \cite{Siebenmann-Sullivan} that the Cannon--Edwards homeomorphism $\Sigma^2 H^3 \to \mathbb S^5$ cannot be taken to be Lipschitz with respect to the polyhedral metric of $\Sigma^2 H^3$, since the suspension circle of $\Sigma^2 H^3$ maps to a set of Hausdorff dimension at least three under any homeomorphism. Since such dimension distortion is possible under quasiconformal mappings, Siebenmann and Sullivan asked whether there still exists a geometric homeomorphism in the form of a quasiconformal mapping $\Sigma^2 H^3 \to \mathbb S^5$. Due to the polyhedral structure, such a  quasiconformal map would be quasisymmetric. 

\section{Absolute continuity of quasisymmetric and quasiconformal mappings}\label{sec:qs}

\begin{question}\label{q:15}
\it If $X$ is a metric space of locally finite Hausdorff $n$-measure and $f$ a quasisymmetric homeomorphism from $X$ onto $\R^n$, $n \ge 2$, is $f$ then absolutely continuous with respect to the Hausdorff $n$-measures?
\end{question}

The answer is NO, as was shown by Romney \cite{rom:singular}.

\begin{question}\label{q:16}
\it If $X$ is a metric space and $f$ a quasisymmetric homeomorphism from $X$ onto $\R^n$, $n \ge 2$, is $f$ then absolutely continuous with respect to the Hausdorff $n$-measures?
\end{question}

The answer is NO. This follows from the ``No'' answer to Question \ref{q:15} given by Romney in \cite{rom:singular}.

\begin{question}\label{q:17}
\it Is there an embedding $f$ of $\R^2$ into $\R^3$ such that, for some $C \ge 1$ and $0<\alpha<1$, $C^{-1} |x-y|^\alpha \le |f(x)-f(y)| \le C|x-y|^\alpha$?
\end{question}

The answer is YES. A construction was given by David and Toro in \cite{dt:reifenberg}.

\begin{question}\label{q:18}
\it Is there a quasisymmetric embedding $f$ of $\R^2$ into $\R^3$ such that the image $f(\R^2)$ contains no rectifiable curves?
\end{question}

The answer is YES. This was first obtained by Bishop in \cite{bis:surface}. It also follows from the `Yes' answer to Question \ref{q:17} given by David and Toro in \cite{dt:reifenberg}.

\subsection{Definitions and notation}\label{subsec:qs-definitions-and-notation}

A measure $\mu$ on a metric space $(X,d)$ is {\it locally finite} if $\mu(B(x,r)) < \infty$ for all balls $B(x,r)$ in $X$. A mapping $f:(X,\mu) \to (Y,\nu)$ between measure spaces is {\it absolutely continuous} if $\mu(E) = 0$ implies $\nu(f(E)) = 0$.

\subsection{Commentary}\label{subsec:qs-commentary}

Questions \ref{q:15} through \ref{q:18} address the measure-theoretic structure of quasiconformal and quasisymmetric images of Euclidean space, especially the problem of absolute continuity of such mappings with respect to the volume measure. This topic has a long history, dating back to foundational work of Gehring and V\"ais\"al\"a in the 1960s.

In \cite{geh:lower-dimensional}, Gehring studied metric and measure-theoretic properties of the images of subspaces under global quasiconformal mappings. In particular, for any quasiconformal map $f:\Omega \to \Omega'$ between domains in $\R^n$, $n \ge 2$, he established the absolute continuity with respect to the Hausdorff $m$-measure $\cH^m$ for the restriction of $f$ to $\Omega \cap V$ for any $m$-dimensional affine subspace $V$, assuming that $\cH^m$ is locally finite on $f(\Omega \cap V)$. V\"ais\"al\"a \cite{vai:embeddings} removed the restriction that the map $f$ arise from a global quasiconformal mapping of an ambient $n$-dimensional domain. More precisely, he showed that any quasisymmetric embedding $f:G \to \R^n$, where $G$ is a domain in $\R^m$ with $2 \le m \le n$ and $\cH^m$ is locally finite on $f(G)$, necessarily is absolutely continuous with respect to $\cH^m$. Finally, in \cite{tys:analytic}, Tyson removed the restriction that the target be contained in any Euclidean space; thus the same conclusion can be stated for any quasisymmetric mapping $f:G \to Y$ onto any metric space such that $\cH^m$ is locally finite on $Y$.

The preceding discussion concerns quasisymmetric mappings $f:G \to Y$ defined on a Euclidean domain $G$ and taking values in a metric space $Y$. Questions \ref{q:15} and \ref{q:16} concern the {\it inverse absolute continuity} problem, where the target is a Euclidean space and the source is a general metric measure space. In this setting it is a priori unclear whether locally finiteness of the relevant Hausdorff measure in the source is needed, hence the formulation in two separate questions.

It is easy to see that quasisymmetric mappings $f:\R^n \to Y$ need not be absolutely continuous with respect to the Hausdorff $n$-measure $\cH^n$ if $\cH^n|_Y$ is not locally finite. A canonical example is the snowflake mapping $\id:\R^n \to (\R^n,d_E^\alpha)$ for $0<\alpha<1$. This map is $\eta$-quasisymmetric with $\eta(t) = t^\alpha$, and the image space $(\R^n,d_E^\alpha)$ has the property that the Hausdorff dimension of any ball $B_{(\R^n,d_E^\alpha)}(x,r)$ is equal to $\tfrac{n}{\alpha}$, whence $\cH^n(B_{(\R^n,d_E^\alpha)}(x,r)) = + \infty$ for all $x \in \R^n$ and $r>0$. 

Assouad's embedding theorem \cite{ass:plongements} implies that any such snowflaked metric space $(\R^n,d_E^\alpha)$ admits a bi-Lipschitz embedding into some finite-dimensional Euclidean space $\R^N$. In general, the ambient dimension $N = N(n,\alpha)$ depends on $n$ and $\alpha$ and can be large. Question \ref{q:17} asks whether there exists a snowflaking parameter $\alpha<1$ for which the space $(\R^2,d_E^\alpha)$ embeds bi-Lipschitzly into $\R^3$, while Question \ref{q:18} asks for a weaker conclusion, namely, an embedding $f$ of $\R^2$ into $\R^3$ with the property that $f(\R^2)$ is free of all nonconstant rectifiable curves. 

We remark in passing here that the precise embedding dimension $N(1,\alpha)$ for bi-Lipschitz embedding of snowflakes of the real line into higher dimensional Euclidean spaces is known; see Assouad \cite[Proposition 4.12]{ass:plongements}. Moreover, there has been substantial interest in estimating the optimal embedding dimension for higher-dimensional Euclidean spaces and for general doubling metric spaces (which is the actual context for Assouad's original theorem). We direct the interested reader to \cite{nn:assouad}, \cite{ds:assouad}, and \cite{tao:assouad} for more information on this latter topic.

A positive answer to Question \ref{q:18} was given by C. Bishop in \cite{bis:surface}. Bishop's construction relies heavily on complex analytic methods, and particularly on the flexibility of two-dimensional conformal mappings. Question \ref{q:17} was subsequently resolved by David and Toro in \cite{dt:reifenberg}. In fact, the main theorem in \cite{dt:reifenberg} is substantially stronger than what is requested in Question \ref{q:17}, as the source space can be any $n$-dimensional {\it Reifenberg flat} metric space. (Without going into details regarding the definition, these are metric spaces that satisfy a bilateral flatness condition similar to the assumption of Reifenberg's Theorem \ref{thm:reifenberg}.) In particular, their theorem covers the case of snowflake embeddings of any Euclidean space $\R^n$, $n \ge 2$, into $\R^{n+1}$, and in \cite[Corollary 2.19 and Remark 2.25]{dt:reifenberg} the authors show the following.

\begin{theorem}[David--Toro]\label{th:david-toro}
Let $n \ge 2$. Then there exists $\alpha = \alpha(n)<1$ and $C>0$ so that for any $R_0>0$, there exists a quasiconformal mapping $f:\R^{n+1} \to \R^{n+1}$ so that $C^{-1} |x-y|^{\alpha} \le |f(x) - f(y)| \le C |x-y|^\alpha$ for all $x,y \in \R^n \times\{0\} \subset \R^{n+1}$ with $|x|,|y| \le R_0$. Moreover, $f$ can be chosen to be equal to the identity in the complement of some neighborhood of $\R^n \times \{0\}$.
\end{theorem}

In his review (\cite{hei:review}) of article \cite{dt:reifenberg} published in Math Reviews as MR1731465, Juha Heinonen wrote: “The well-known von Koch snowflake, a self-similar Jordan curve in the plane, is a prototypical fractal and a source of numerous examples and counterexamples in geometric analysis. By following its iterative construction, it is easy to produce a parameterization of the snowflake curve by the round circle via a homeomorphism$f$ satisfying $|f(x)-f(y)| \approx |x-y|^\alpha$, where $0<\alpha<1$ is the reciprocal of the Hausdorff dimension of the curve \dots It has been a long-standing open problem how to construct analogous ‘snowballs’ in higher dimensions. The wonderful paper under review accomplishes this and more."

Among many applications and consequences, we note here that the examples constructed by David and Toro were used by Kaufman, Tyson, and Wu \cite{ktw:qr} in the construction of $C^{1,\alpha}$ branched quasiregular mappings of $\R^n$, $n \ge 4$, extending an earlier three-dimensional construction of Bonk and Heinonen \cite{bh:qr} and answering a question posed by V\"ais\"al\"a \cite{vai:icm} in his 1978 ICM address. For more information about quasiregular mappings and their branch sets, see section \ref{sec:qr}.

The full inverse absolute continuity problem remained open for many years prior to its solution by Romney. We observe in particular that the authors of \cite{abh:sobolev-bv} note in the introduction that the results in that paper were strongly motivated by efforts to solve the inverse absolute continuity problem. The eventual solution to the problem resides in the following two theorems of Matthew Romney \cite{rom:singular}.\footnote{Kai Rajala and Matthew Romney also shared an interesting anecdote with the authors. Rajala's work \cite{raj:spheres} on geometrically quasiconformal parameterizations of spheres (discussed above in connection with Theorem \ref{thm:bk}) was inspired by Gehring's inverse absolute continuity problem, while Romney's solution to that problem in \cite{rom:singular} in turn drew inspiration from the results in \cite{raj:spheres}.}

\begin{theorem}[Romney]\label{th:romney-1}
For each $n \ge 2$ there exists a metric space $(X,d)$, a quasisymmetric homeomorphism $f:[0,1]^n \to X$, and a set $E \subset X$ so that $\cH^n(X \setminus E) = 0$ and $\cH^n(f(E)) = 0$.
\end{theorem}

\begin{theorem}[Romney]\label{th:romney-2}
For each $n \ge 2$ and $0<\delta<1$ there exists a metric space $(X,d)$ with $\cH^n(X)<\infty$, a quasisymmetric homeomorphism $f:[0,1]^n \to X$, and a set $E \subset X$ so that $\cH^n(E) \ge \delta$ and $\cH^n(f(E)) = 0$.
\end{theorem}

An interesting consequence of Theorem \ref{th:romney-2}, which answers \cite[Question 17.3]{raj:spheres} is the following observation: {\it There exists a metric space $(X,d)$ so that $\cH^2|_X$ is locally finite, and a quasisymmetric mapping $f:[0,1]^2 \to X$ so that $f$ is not geometrically quasiconformal.}

To conclude this section, returning full circle to Gehring's original work, we mention subsequent work of Ntalampekos and Romney in \cite{nr:inverse} which gave a negative resolution to the inverse absolute continuity problem in its original setting, namely, the problem of absolute continuity of globally defined quasiconformal maps of $\R^3$ with respect to the Hausdorff $2$-measure on hyperplanes. More precisely, in \cite{nr:inverse} the authors construct a quasiconformal map $f:\R^3 \to \R^3$ and a Borel set $E \subset \R^{2} \times \{0\}$ so that $\cH^2(E)>0$ and $\cH^2(f(E)) = 0$.

\section{Analysis on metric measure spaces}\label{sec:mms}

\begin{question}\label{q:19}
\it If $X$ is an Ahlfors $Q$-regular space that admits a weak $(1,1)$-Poincar\'e inequality, is $Q$ then an integer?
\end{question}

The answer is NO, by results of Bourdon and Pajot \cite{bp:buildings} and Laakso \cite{laa:PI}.

\begin{question}\label{q:20}
\it If $X$ is an Ahlfors $Q$-regular Loewner space for some $Q>1$, is $Q$ then an integer?
\end{question}

The answer is NO. This follows from the answer to Question \ref{q:19} and results in \cite{hk:quasi}.

\subsection{Definitions and notation}\label{subsec:mms-definitions-and-notation}

 A nonnegative Borel function $g$ on a metric space $(X,d)$ is said to be an {\it upper gradient} of a function $u:(X,d) \to \R$ if the inequality
$$
|u(x) - u(y)| \le \int_\gamma g \, ds
$$
holds for all points $x,y \in X$ and all rectifiable curves $\gamma$ in $X$ such that $\gamma$ joins $x$ to $y$. Here $\int_\gamma g \, ds$ denotes the line integral of $g$ along $\gamma$, computed with respect to the arc length parameterization of $\gamma$.

An Ahlfors $Q$-regular metric measure space $(X,d,\mu)$ is said to satisfy the {\it weak $(1,1)$-Poincar\'e inequality} if there exist constants $C>0$ and $\tau \ge 1$ so that the inequality
$$
\int_{B(x,r)} |u-u_{B(x,r)}| \, d\mu \le C r \int_{B(x,\tau r)} g \, d\mu
$$
holds true for all balls $B(x,r) \subset X$ with $0<r<\diam(X)$, all continuous functions $u:B(x,\tau r) \to \R$, and all upper gradients $g$ of $u$ in $B(x,\tau r)$. 

A {\it ring domain} in $(X,d)$ is a pair $(E,F)$, where $E$ and $F$ are disjoint, nondegenerate continua in $X$, and the {\it conformal modulus} of a ring domain $(E,F)$ is
$$
\Modd_Q(E,F) := \inf \int_X g^Q \, d\mu
$$
where the infimum is taken over all nonnegative Borel functions $g$ so that $\int_\gamma g \, ds \ge 1$ for all rectifiable curves $\gamma$ such that $\gamma$ joins a point of $E$ to a point of $F$. With this notion in place, an Ahlfors $Q$-regular metric measure space $(X,d,\mu)$ is said to be {\it Loewner} if the function
$$
\lambda_X(t) := \inf \{ \Modd_Q(E,F): \dist(E,F) \le t \min \{ \diam(E),\diam(F) \}
$$
is positive for each $t>0$.

\subsection{Commentary}\label{subsec:mms-commentary}

Questions \ref{q:19} and \ref{q:20} were among the earliest problems in the list to be resolved. At the time of publication of \cite{hs:questions}, the class of Ahlfors regular spaces satisfying the weak Poincar\'e inequality, or alternatively the Loewner condition, was already known to be relatively broad. A list of examples in \cite[Section 6]{hk:quasi} includes Euclidean space and compact Riemannian manifolds, the Heisenberg group and other sub-Riemannian Carnot groups, noncompact Riemannian manifolds with nonnegative Ricci curvature, complete, connected, orientable, and linearly locally contractible topological manifolds, and various simplicial complexes. In addition, Section 6 of \cite{hk:quasi} describes various deformation and gluing procedures which preserve either the weak Poincar\'e inequality or the Loewner condition. A common feature of all of these examples was that the Hausdorff dimension was an integer, and the various deformation and gluing procedures described above preserved the integrality of Hausdorff dimension. This observation motivated the statements of Questions \ref{q:19} and \ref{q:20}, which essentially ask in what breadth the newly developed theory of analysis in metric spaces applies.

Following Mostow's rigidity theorem \cite{mos:rigidity1}, \cite{mos:rigidity2} and subsequent work by Kor\'anyi--Reimann, Pansu, and Bourdon, the role of abstract metric analysis on the boundaries of Gromov hyperbolic groups and spaces for the large-scale geometric properties of the underlying hyperbolic space was well understood. Various constructions of Fuchsian hyperbolic buildings had already been considered by Bourdon \cite{bou1}, \cite{bou2}. In particular, in \cite{bou2} Bourdon defined a family of Fuchsian hyperbolic buildings $I_{p,q}$, indexed by a pair of integers $p \ge 5$ and $q \ge 3$. The space $I_{p,q}$ is obtained by gluing together right-angled regular hyperbolic $p$-gons along edges, subject to the constraint that the link at each vertex is the complete bipartite graph $K_{q,q}$. In \cite{bou2} Bourdon showed that the space $I_{p,q}$ satisfies the $CAT(-1)$ condition and consequently is Gromov hyperbolic. He also established the Ahlfors $Q_{p,q}$-regularity of the boundary $\partial_\infty I_{p,q}$ for a specific choice of visual metric, and observed that this metric infimized dimension among all such visual metrics. (In modern language, the latter result states that Bourdon's metric realizes the {\em Ahlfors regular conformal dimension} of the boundary.) The stage was therefore naturally set for Bourdon and Pajot's solution to Question \ref{q:19} in \cite{bp:buildings}.

\begin{theorem}[Bourdon--Pajot]
Let $p \ge 5$ and $q \ge 3$ be integers, and let $I_{p,q}$ denote the Fuchsian hyperbolic building constructed in \cite{bou2}. Let $\partial_\infty I_{p,q}$ denote the Gromov boundary of $I_{p,q}$, equipped with the metric $d$ defined in \cite{bou2}, and let $\mu$ denote the Hausdorff $Q_{p,q}$-measure on $\partial_\infty I_{p,q}$. Then $(\partial_\infty I_{p,q},d,\mu)$ is an Ahlfors regular space supporting the weak $(1,1)$-Poincar\'e inequality. Moreover, $(\partial_\infty I_{p,q},d,\mu)$ also satisfies the Loewner condition.
\end{theorem}

The value of the dimension of $\partial_\infty I_{p,q}$ was already computed in \cite{bou2}; this value is
$$
Q_{p,q} = 1 + \frac{\log(q-1)}{\arccosh(\tfrac{p-2}{2})}.
$$
It is a straightforward number-theoretic exercise to show that $Q_{p,q}$ is non-integer (in fact, irrational) for all relevant $p$ and $q$. As observed in \cite[Example 11.2]{tys:thesis}, the set $\{Q_{p,q}:p \ge 5, q \ge 3\}$ is dense in the interval $[1,\infty)$, and $Q_{p,q} = Q_{p',q'}$ if $q'=1+(q-1)^N$ and $p'=2+2T_N(p/2-1)$ for some integer $N$, where $T_N$ denotes the $N$th Chebyshev polynomial. For example, $Q_{5,3} = Q_{9,5}$.

Subsequently, Laakso provided examples of Ahlfors $Q$-regular spaces supporting the weak $(1,1)$-Poincar\'e inequality for arbitrary real numbers $Q \in (1,\infty)$. In \cite{laa:PI}, Laakso considers metric spaces obtained as quotients of product spaces $Z \times [0,1]$ for suitable spaces $Z$ of Cantor type. The relevant quotient operation consists in identifying suitable pairs $(z_1,t_1),(z_2,t_2) \in Z \times [0,1]$, and defining a quotient metric $d$ on $X = Z \times [0,1] / \sim$ as follows: for two equivalence classes $[(z_i,t_i)]$ and $[(z_f,t_f)]$ in $X$, the distance $d([(z_i,t_i)],[(z_f,t_f)])$ is the infimum of the total length of all countable unions $A$ of line segments in $Z \times [0,1]$ such that (a) $[(z_i,t_i)]$ and $[(z_f,t_f)]$ are contained in $A$, and (b) the final endpoint of the $k$th segment in $A$ is identified with the initial endpoint of the $(k+1)$st segment. Laakso calls the identifications in $Z \times [0,1]$ {\em wormholes}. The main theorem of \cite{laa:PI} reads as follows:

\begin{theorem}[Laakso]
Let $1<Q<\infty$. There exists a wormhole construction as above so that the quotient metric space $(X,d)$ is Ahlfors $Q$-regular and supports the weak $(1,1)$-Poincar\'e inequality, when equipped with the Hausdorff measure $\cH^Q$.
\end{theorem}

As a side note, we comment that Laakso's proof that $(X,d)$ is Ahlfors $Q$-regular relies on the fact that the quotient map $Z \times [0,1] \to X$ is David--Semmes regular; see Question \ref{q:8} for the definition. Laakso also proves (Theorem 4.1 in \cite{laa:PI}) that none of his wormhole spaces admit any bi-Lipschitz embedding into any finite-dimensional Euclidean space.

\section{Analysis and geometry in the Heisenberg group}\label{sec:heisenberg}

\begin{question}\label{q:21}
\it If the first Heisenberg group is deformed by a metric doubling measure, is a $(1,1)$-Poincar\'e inequality retained?
\end{question}

This question remains open.

\begin{question}\label{q:22}
\it Is the first Heisenberg group minimal in looking down?
\end{question}

The answer is NO, by a result of Le Donne, Li, and T. Rajala \cite{llr:heisenberg}.

\begin{question}\label{q:23}
\it If two $Q$-regular BPI metric spaces are look-down equivalent, are they then BPI equivalent?
\end{question}

According to an unpublished preprint of Laakso \cite{laa:lookdown-bpi}, the answer is NO.

\begin{question}\label{q:24}
\it If $f$ is a Lipschitz map of a subset $E$ in the first Heisenberg group equipped with its Carnot metric into a metric space such that the Hausdorff $4$-measure of $f(E)$ is positive, is $f$ then bi-Lipschitz in a subset of $E$ of positive Hausdorff $4$-measure?
\end{question}

The answer is NO. This follows from the `No' answer to Question \ref{q:22} given by Le Donne, Li, and Rajala in \cite{llr:heisenberg}. A more direct construction which answers Question \ref{q:24} directly can also be found in that paper.

\begin{question}\label{q:25}
\it Can the $t$-axis in the first Heisenberg group be mapped onto a locally rectifiable curve by a quasiconformal self-map of the group?
\end{question}

This question remains open.

\subsection{Definitions and notation}\label{subsec:heisenberg-definitions-and-notation}

The {\it first Heisenberg group} $\Heis^1$ is the unique connected and simply connected three-dimensional nonabelian stratified Lie group. For the purposes of this survey, we present $\Heis^1$ as the Euclidean space $\R^3 = \C \times \R$ with coordinates $\bp = (z,t) = (x,y,t)$ and group law $(z,t) * (z',t') = (z+z',t+t'+2\Imag(z \overline{z'}))$. The Heisenberg group is equipped with a sub-Riemannian structure coming from the two-dimensional, left invariant nonintegrable distribution $H\Heis^1$ whose fiber at a point $\bp \in \Heis^1$ is $H_\bp \Heis^1 = \spa \{ X(\bp), Y(\bp) \}$, where $X = \partial_x + 2y \partial_t$ and $Y = \partial_y - 2x \partial t$. An absolutely continuous curve $\gamma:I = (a,b) \to \Heis^1$ is said to be {\it horizontal} if $\gamma'(s) \in H_{\gamma(s)} \Heis^1$ for a.e.\ $s \in I$; the Chow--Rashevsky theorem guarantees that $\Heis^1$ is {\it horizontally connected}, i.e., every pair of points $\bp,\bq \in \Heis^1$ can be joined by at least one horizontal curve. (In fact, it is relatively easy to construct such a curve as a concatenation of several integral curves associated to the vector fields $X$ and $Y$.) Defining the {\it horizontal length} of a horizontal curve $\gamma$ to be $\ell_h(\gamma) := \int_a^b \sqrt{c(s)^2 + d(s)^2} \, ds$, where $\gamma'(s) = c(s) X_{\gamma(s)} + d(s) Y_{\gamma(s)}$ a.e., we introduce the {\it Carnot--Carath\'eodory (C--C) metric} $d_{cc}(\bp,\bq) := \inf \ell_h(\gamma)$ where the infimum is taken over all horizontal curves $\gamma$ joining $\bp$ to $\bq$. 

The metric space $(\Heis^1,d_{cc})$ admits a one-parameter family of scaling automorphisms $(\delta_r)_{r>0}$ where $\delta_r:\Heis^1 \to \Heis^1$ is given by $\delta_r(z,t) = (rz,r^2t)$. Noting that the Haar measure of $\Heis^1$ coincides (up to a constant multiple) with the Lebesgue measure on the underlying Euclidean space, and computing the Jacobian determinant of the automorphism $\delta_r$, $r>0$, one quickly deduces that the Hausdorff dimension of $(\Heis^1,d_{cc})$ is equal to $4$. The metric measure space $(\Heis^1,d_{cc},\cH^4_{d_{cc}})$ is Ahlfors $4$-regular and satisfies the $(1,1)$-Poincar\'e inequality. 

The $t$-axis $\{ \bp = (0,t) \in \Heis^1 : t \in \R \}$ is a nonrectifiable curve in $(\Heis^1,d_{cc})$. In fact, there exists a constant $c>0$ so that $d_{cc}((0,t_1),(0,t_2)) = c \sqrt{|t_2-t_1|}$ for all $t_1,t_2 \in \R$. It follows that the Hausdorff dimension of the $t$-axis in the C--C metric is equal to two.

For the definition of metric doubling measure, see subsection \ref{subsec:qc-uniformization-parametrization-definitions-and-notation}.

A map $g:(X,d) \to (X',d')$ between metric spaces is said to be {\it $C$-conformally bi-Lipschitz with scale factor $\lambda>0$} if $g$ is $C$-bi-Lipschitz as a map from the rescaled metric space $(X,\lambda d)$ to $(X',d')$. That is, for all $x,y \in X$, we have $C^{-1} \lambda d(x,y) \le d'(f(x),f(y)) \le C \lambda d(x,y)$.

An Ahlfors $Q$-regular metric space $(X,d)$ is said to be a {\it BPI space}\footnote{The acronym BPI stands for `big pieces of itself'.} if there exist constants $C \ge 1$ and $c>0$ so that given any two balls $B(x_1,r_1)$ and $B(x_2,r_2)$ in $X$ with $0<r_1,r_2<\diam(X)$, there is a closed subset $A$ of $B(x_1,r_1)$ with the property that $\cH^Q(B(x_1,r_1) \ge c r_1^Q$ and $A$ is $C$-conformally bi-Lipschitz to a subset of $B(x_2,r_2)$ with scale factor $\tfrac{r_2}{r_1}$. Two BPI spaces $(X,d)$ and $(X',d')$ of the same regularity dimension $Q$ are said to be {\it BPI equivalent} if there are constants $C\ge 1$ and $c>0$ so that given any two balls $B(x,r) \subset X$ and $B(x',r') \subset X'$, there is a closed subset $A$ of $B(x,r)$ with the property that $\cH^Q(B(x,r)) \ge c r^Q$ and $A$ is $C$-conformally bi-Lipschitz to a subset of $B(x',r')$ with scale factor $\tfrac{r'}{r}$.

Given two BPI spaces $X$ and $Y$, both with regularity dimension $Q$, we say that $X$ {\it looks down} on $Y$ if there exists a closed set $A \subset X$ and a Lipschitz map $f:A \to Y$ so that $\cH^Q_{d'}(f(A))>0$. If $X$ looks down on $Y$ and $Y$ also looks down on $X$< we say that $X$ and $Y$ are {\it look-down equivalent}. Finally, a BPI space $X$ of regularity dimension $Q$ is said to be {\it minimal in looking down} if for every BPI space $Y$ of regularity dimension $Q$ so that $X$ looks down to $Y$, it holds true that $X$ is BPI equivalent to $Y$.

Quasiconformality and quasisymmetry were defined in subsection \ref{subsec:qc-uniformization-parametrization-definitions-and-notation} for homeomorphisms of Euclidean domains, and these metric notions extend {\it mutatis mutandis} to the Heisenberg group equipped with its Carnot--Carath\'eodory metric. There is also an analytic definition for quasiconformality in this context; as this notion would take time to define and is not relevant for the discussion at hand we omit it. As in the Euclidean case, the analytic and metric definitions of quasiconformality are equivalent, and are in turn also equivalent to local quasisymmetry.

\subsection{Commentary}\label{subsec:heisenberg-commentary}

Questions \ref{q:22}, \ref{q:23}, and \ref{q:24} are closely related to  Question \ref{q:9} and the ensuing discussion in subsection \ref{subsec:bl-parametrizations-embeddings-commentary}. Schul's positive answer to Question \ref{q:9} relied on a `quantitative inverse function theorem' for Lipschitz maps from Euclidean spaces to arbitrary metric spaces, stated above as Theorem \ref{thm:schul}. If a similar result held for Lipschitz mappings from the Heisenberg group or between BPI spaces, then Questions \ref{q:22} - \ref{q:24} would have positive answers. It turns out, however, that it does not.

In \cite{llr:heisenberg}, Le Donne, Li, and Rajala construct a Lipschitz mapping (in fact, a Lipschitz homeomorphism) from the Heisenberg group onto another Ahlfors $4$-regular space that fails to be bi-Lipschitz on \emph{any} positive-measure subset of the Heisenberg group. This immediately shows that even a qualitative version of Theorem \ref{thm:schul} fails for mappings from the Heisenberg group, and that Question \ref{q:24} has a negative answer.

The key property of the Heisenberg group with its Carnot--Carath\'eodory metric $d_{cc}$ (which Euclidean space lacks) that allows for their construction is the following ``shortcut property'': There is a constant $\lambda\in (0,1)$ such that for all $p\in \Heis^1$ and all $r>0$, there are points $q_1, q_2 \in B(p,r)$ such that
$$ d_{cc}(q_1,q_2) \geq \lambda r$$
and
$$ d_{cc}(p_1, p_2) \leq d_{cc}(p_1, q_1) + d_{cc}(p_2, q_2) \text{ for all } p_1, p_2 \notin B(p,r).$$
In other words, the points $q_1$ and $q_2$ allow for a ``shortcut'' from $p_1$ to $p_2$ through $B(p,r)$ that cuts off a (substantial) term of the triangle inequality. The choice of $q_1$ and $q_2$ is simple. If $p=(0,0,0)$ and $r=1$, one can take $q_1 = p$ and $q_2 = (0,0,1/4)$, and for other balls one can reduce to this case by scaling and left translation.

By cleverly exploiting these shortcuts in $\Heis^1$ at all locations and scales, Le Donne, Li, and Rajala are able to construct a new Ahlfors $4$-regular metric $d$ on the Heisenberg group that is bounded above by $d_{cc}$ but is not bi-Lipschitz equivalent to it on any positive-measure set, answering Question \ref{q:24}. A more careful, self-similar version of their construction allows the metric $d$ to be BPI, like $d_{cc}$ itself, and therefore also answers Question \ref{q:22} in the negative.

Question \ref{q:23} is answered in a strong sense by Laakso in \cite{laa:lookdown-bpi}, which has unfortunately never been formally published. By using two similar variations of his `wormhole' construction (see subsection \ref{subsec:mms-commentary}), Laakso constructs two BPI spaces that are look-down equivalent and even support $(1,1)$-Poincar\'e inequalities, but that fail to be BPI equivalent. Recently it has been fruitful to recast and generalize some of Laakso's constructions, e.g., those in \cite{laa:PI}, as inverse limits of metric graphs. (See for instance \cite{ck:graphs} and \cite{aeb:graphs}.) It would be interesting to rework some of these mapping theory examples in that context as well.

This group of questions seems to admit a natural follow up. We know from the discussion in subsection \ref{subsec:bl-parametrizations-embeddings-commentary} that Lipschitz mappings from $\mathbb{R}^n$ to arbitrary metric spaces admit both qualitative and quantitative metric inverse function theorems, of the type that say that Lipschitz mappings admit Lipschitz inverses on positive-measure fractions of their image. Such results therefore also hold for domains that are mild generalizations of Euclidean space, like rectifiable or uniformly rectifiable spaces. On the other hand, we also know that such theorems may fail when the domain of the mapping is otherwise quite reasonable, like the Ahlfors regular, BPI, PI space examples in Questions \ref{q:23} and \ref{q:24}. We might therefore ask: \emph{does this type of (qualitative or quantitative) metric inverse function theorem characterize rectifiable or uniformly rectifiable metric spaces?}\footnote{Just as this survey was in its final stages of preparation, Li and Schul posted \cite{lischul}, answering the qualitative version of this question.}

Question \ref{q:25} may strike the reader as somewhat specialized, however, at its core this question concerns the flexibility of Heisenberg quasiconformal mappings. One of the primary difficulties which stymied work on this question is the comparative lack (relative to the Euclidean or Riemannian case) of methods for constructing explicit quasiconformal maps in the Heisenberg group. We briefly discuss some standard methods for constructing such maps, and the extent to which such methods do or do not remain available in the sub-Riemannian Heisenberg setting.

Reimann \cite{rei:flow} showed that the $s$-advance homeomorphism $p \mapsto \gamma_p(s)$, associated to an ODE in $\R^n$ of the form $\dot\gamma_p(s) = V(\gamma_p(s))$, $\gamma_p(0) = p$, is quasiconformal provided the vector field $V$ has finite {\it Reimann $Q$-seminorm}. Reimann's seminorm $||V||_Q$, whose definition we omit, lies between the Lipschitz seminorm $||V||_L := \sup_{x \ne x'} |V(x)-V(x')|/|x-x'|$ (the finiteness of which implies that the advance homeomorphisms are bi-Lipschitz) and the Zygmund seminorm $||V||_Z := \sup_x \sup_{y \ne 0} |V(x+y)+V(x-y)-2V(x)|/|y|$ (the finiteness of which implies that the advance homeomorphsims are bi-H\"older). McMullen \cite[Appendix A]{mcm:3manifolds} provides a succinct description of Reimann's flow method for the construction of (Euclidean) quasiconformal mappings. 

An analogous methodology exists in the Heisenberg setting, as shown by Kor\'anyi and Reimann in \cite{kr:foundations}. This allows for the construction of certain explicit smooth quasiconformal homeomorphisms of Heisenberg domains. However, it is not known whether every quasiconformal map in $\Heis^1$ may be embedded in such a flow.\footnote{Note that a similar state of affairs holds true in $\R^n$ for $n \ge 3$. In $\R^2$, solvability of the Beltrami equation $f_{\overline{z}} = \mu f_z$ via the measurable Riemann mapping theorem implies that every quasiconformal map of planar domains may be embedded in a flow of quasiconformal maps.}

The fact that quasiconformal mappings are locally quasisymmetric permits the use of various convergence theorems, allowing for the construction of nonsmooth quasiconformal maps as locally uniform limits of sequences of maps with a uniform bound on quasiconformal dilatation and suitable normalization conditions. Such convergence results also hold true in the Heisenberg setting (and in more general metric measure spaces).

What remains unknown in the latter setting, however, is any analog for the natural PL constructions which feature prominently in numerous explicit examples of quasiconformal and quasiregular maps. Inherent in many of these constructions is the fact that the reflection map $(x_1,x_2,\ldots,x_n) \mapsto (-x_1,x_2,\ldots,x_n)$ is a quasiconformal self-map of $\R^n$. But the anisotropic character of the Carnot--Carath\'eodory metric, and in particular the role of the contact geometry of the Heisenberg group, precludes similar constructions.

Using the Kor\'anyi--Reimann flow method and compactness arguments, Balogh \cite{bal:hausdorff} constructed, for each $0<\alpha<1$ and each $0<\beta<2$, a subset $A$ of a rectifiable curve in $\Heis^1$ with Hausdorff dimension $\alpha$, and a quasiconformal self-map $f:\Heis^1 \to \Heis^1$ so that $f(A)$ lies in the $t$-axis and the Hausdorff dimension of $f(A)$ is $\beta$. While this result shows that large subsets of rectifiable curves can be quasiconformally mapped onto large subsets of the $t$-axis, it appears unlikely that the techniques of the proof can fully answer Question \ref{q:25}.

\section{Quasiconformal and quasiregular mappings on smooth manifolds 
}\label{sec:qr}

\begin{question}\label{q:26}
\it Is the fixed point set of every quasiconformal reflection $f:\Sph^3 \to \Sph^3$ a topologically tame $2$-sphere?
\end{question}

This question remains open.

\begin{question}\label{q:27}
\it Can every branched cover $\Sph^n \to \Sph^n$, $n \ge 3$, be made BLD by changing the metric in the domain but keeping the space $n$-regular and linearly locally contractible?
\end{question}

This question remains open.

\begin{question}\label{q:28}
\it Is every branched cover $f:\Sph^n \to \Sph^n$, $n \ge 3$, topologically conjugate to a quasiregular map?
\end{question}

This question remains open.

\begin{question}\label{q:29}
\it Is the topological dimension of the branch set of a quasiregular map $\Sph^n \to \Sph^n$, $n \ge 3$, either $-1$, $n-2$, or $n$?
\end{question}

This question remains open.

\begin{question}\label{q:30}
\it Does there exist a nonconstant quasiregular map $\R^4 \to \Sph^2 \times \Sph^2 \# \Sph^2 \times \Sph^2$?
\end{question}

The answer is YES. This was established by Rickman in \cite{rick:qr-elliptic-4-manifolds}.

\subsection{Definitions and notation}\label{subsec:qr-definitions-and-notation}

A continuous mapping $f\colon \mathbb S^3 \to \mathbb S^3$ is an \emph{involution} if $f \circ f = \id$. Since $f^{-1} = f$, an involution is a homeomorphism. An orientation-reversing involution is called a \emph{reflection}; in other words, a reflection is an orientation-preserving involution with respect to opposite orientations. 

A fixed point set $\mathrm{Fix}(f)$ of a self-map $f\colon \Sph^3 \to \Sph^3$ is the set $\mathrm{Fix}(f) = \{ x\in \Sph^3 \colon f(x)=x\}$.  
An embedded $k$-sphere $S\subset \Sph^3$, for $k\in \{1,2\}$, is \emph{tame} if there exists a homeomorphism $\varphi \colon \Sph^3 \to \Sph^3$ for which $\varphi(S) = \Sph^k \subset \Sph^3$; an embedded $k$-sphere in $\Sph^3$ is \emph{wild} if it is not tame. Furthermore, an involution $f \colon \Sph^3 \to \Sph^3$ is \emph{tame} if its fixed point set $\mathrm{Fix}(f)$ is a tame $k$-sphere for $k\in \{1,2\}$; an involution $\Sph^3 \to \Sph^3$ is \emph{wild} if it is not tame.

A mapping $f\colon \Sph^n\to \Sph^n$ is \emph{discrete} if the fibers $f^{-1}(y)$ are discrete sets for $y\in Y$, and \emph{open} if the image $fU$ of an open set $U$ is open. A continuous mapping $f\colon \Sph^n \to \Sph^n$ is a \emph{branched cover} if it is an orientation-preserving discrete and open mapping. 

A branched cover $f\colon \Sph^n \to \Sph^n$ has \emph{bounded length distortion (BLD)} if there exists a constant $L\ge 1$ having the property that, for every curve $\gamma$ in $\Sph^n$, $\frac{1}{L} \ell(\gamma) \le \ell(f\gamma) \le L \ell(\gamma)$, where $\ell(\cdot)$ is the length of a curve.

A mapping $f\colon M \to N$ between oriented Riemannian $n$-manifolds is \emph{quasiregular} if $f$ belongs to the Sobolev space $W^{1,n}_{\mathrm{loc}}(M,N)$ and there exists a constant $K\ge 1$ for which $\norm{Df}^n \le K J_f$ almost everywhere in $M$, where $\norm{Df}$ is the operator norm and $J_f = \star f^*\mathrm{vol}_N$ the Jacobian determinant of the differential $Df$. Mappings $f_1\colon \Sph^n \to \Sph^n$ and $f_2 \colon \Sph^n \to \Sph^n$ are \emph{conjugate} if there exists homeomorphisms $\phi \colon \Sph^n \to \Sph^n$ and $\psi \colon \Sph^n \to \Sph^n$ for which $f_2 = \psi \circ f_1 \circ \phi$.

A point $x\in M$ is a \emph{branch point of a quasiregular map $f\colon M \to N$} if $f$ is not a local homeomorphism at $x$. The \emph{branch set $B_f$ of $f$} is the set of all branch points of $f$.

\subsection{Commentary}\label{subsec:qr-commentary}

To set the scene for Question \ref{q:26}, it is a classical theorem of P.A.~Smith \cite{Smith} that the fixed point set $\mathrm{Fix}(f)$ of an involution $f\colon \Sph^3 \to \Sph^3$ is a topological $2$-sphere if $f$ is orientation reversing and a topological $1$-sphere if $f$ is orientation preserving. By a celebrated theorem of Bing \cite{bing:double}, which originally motivated the introduction of the decomposition space $\R^3/\cB$ as in Question \ref{q:11}, there exists an involution $\Sph^3 \to \Sph^3$ whose fixed point set is an Alexander horned sphere, a famous wildly embedded $2$-sphere in $\Sph^3$ due to Alexander \cite{Alexander} in 1920.

The picture, however, changes under for smooth involutions. By Bochner's theorem \cite{Bochner}, the fixed point set of a $C^1$-smooth involution is a smooth submanifold, and hence locally tame. Thus, by Morton Brown's solution to the generalized Sch\"oenflies conjecture \cite{Brown}, $C^1$-smooth reflections are tame. Finally, by solution of the Smith conjecture, all $C^1$-smooth involutions $\Sph^3 \to \Sph^3$ are tame; see the festschrift of Morgan and Bass \cite{Morgan-Bass}. 

Question \ref{q:26} may now be seen as a question whether the Smith conjecture holds in the quasiconformal category. In two-dimensions the problem is completely understood. Indeed, from topological point of view, the fixed point set of a reflection $\Sph^2 \to \Sph^2$ is a Jordan curve and hence tame by Sch\"oenflies theorem. Regarding quasiconformal reflections, Ahlfors gave in \cite{Ahlfors-Acta-1963} -- based on previous work with Beurling \cite{Beurling-Ahlfors} -- a characterization of Jordan curves, which are fixed point sets of quasiconformal reflections $\Sph^2 \to \Sph^2$. 

Coming back to Question \ref{q:26}, Hamilton announced in \cite{Hamilton-Contemp-Math} the following solution to this problem: \emph{the complementary domains of a fixed point set $\mathrm{Fix}(f)$ of a quasiconformal reflection $\Sph^3 \to \Sph^3$ are simply-connected}. We are not aware that a complete proof of this solution has appeared in print. To our knowledge, the best published result in this positive direction is a theorem of Grillet \cite{Grillet}, which implies that \emph{a bi-Lipschitz reflection $\Sph^3 \to \Sph^3$ with small bi-Lipschitz constant is tame}. Grillet's theorem reads as follows.

\begin{theorem}[Grillet]
Any action of a finite cyclic group by $(1+\frac{1}{4000})$-bi-Lipschitz homeomorphisms on a closed $3$-manifold is conjugate to a smooth action. 
\end{theorem}

Regarding the quasiconformal tameness, it is known, by a result of Heinonen and Yang \cite{Heinonen-Yang}, that complementary components of the fixed point set of a quasiconformal reflection are uniform domains. 

Question \ref{q:26} may also be seen in particular as a question about the regularity of Bing's wild involution, i.e., whether Bing's involution has a quasiconformal representative. In this direction, Freedman and Starbird in \cite{Freedman-Starbird}\footnote{Unpublished at the time of writing.} discuss gauge of modulus of continuity functions related to Bing's construction and show that the gauge does not allow H\"older continuous reflections. This complements Gehring's classical result that a quasiconformal map $\Sph^3 \to \Sph^3$ is H\"older continuous, and yields a partial positive solution to Question \ref{q:26}: \emph{Bing's involution does not have a quasiconformal representative}. We note in passing that it is, however, possible to construct Bing's involutions $\Sph^3 \to \Sph^3$ in Sobolev spaces $W^{1,p}(\Sph^3,\Sph^3)$ for $p<2$; see \cite{Onninen-Pankka} for discussion.  

Questions \ref{q:27}, \ref{q:28}, and \ref{q:29} are tied together by a classical theorem of Reshetnyak \cite{Reshetnyak}: \emph{every non-constant quasiregular mapping is a branched cover}. Apart from Question \ref{q:28} we are not aware of progress in these questions and the ICM survey of Heinonen \cite{Heinonen-ICM} from 2002 is still mostly up to date. We note in passing that the term {\it bounded length distortion} in Question \ref{q:27} was coined originally by Martio and V\"ais\"al\"a in \cite{Martio-Vaisala} for a subclass of quasiregular mappings, but its importance in the metric context was understood after works of Heinonen and Rickman \cite{Heinonen-Rickman}, \cite{Heinonen-Rickman_Duke} and the work of Heinonen and Sullivan \cite{Heinonen-Sullivan} on branched parametrizations of metric spaces. See e.g.~Drasin--Pankka \cite{Drasin-Pankka}, Heinonen--Keith \cite{Heinonen-Keith}, or Onninen--Rajala \cite{Onninen-Rajala} for other applications. 

Questions \ref{q:27} and \ref{q:28} may be seen as uniformization problems for general topological branched covers. Question \ref{q:27} bear resemblance to parametrization questions in geometric topology, Questions \ref{q:10}--\ref{q:15}, as they ask whether general branched covers have good internal geometry. For $n=2$ such good internal geometry exists by classical Sto\"ilow's theorem  \cite{Stoilow-1928,Stoilow-1956}: \emph{a branched cover $f \colon S \to S'$ of Riemann surfaces is a composition $f = h\circ \varphi$ of a holomorphic $h \colon S\to S'$ and a quasiconformal mapping $\varphi \colon S\to S$}; see e.g.~the monograph of Astala, Iwaniec, and Martin \cite{AIM-book}. In particular, the branch set in this case is always discrete.

 


Note that the branch set of a quasiregular mapping need not be contained in a simplicial complex. Indeed, in \cite{Heinonen-Rickman} Heinonen and Rickman construct an elaborate example of a quasiregular map $\Sph^3 \to \Sph^3$ whose branch set consists of an Antoine's necklace and an infinite $1$-complex tending to the necklace. This was also a first example of a quasiregular mapping whose Jacobian is not (locally) comparable to the Jacobian of a quasiconformal mapping.

Question \ref{q:29} stems from a search for a higher dimensional counterpart of the topological version of Sto\"ilow's theorem: \emph{light open mappings between surfaces are discrete mappings with discrete branch sets}; see also \cite{Luisto-Pankka} for a discussion on the proof. Recall that a mapping $f$ is \emph{light} if its pre-image fibers $f^{-1}(y)$ are totally disconnected. For $n\ge 3$, the topological counterpart of Stoilow's theorem is the classical {\v C}ernavski{\u \i}--V\"ais\"al\"a theorem (\cite{Cernavskii, Vaisala-1966}), which states that \emph{the branch set of a branched cover between manifolds has codimension at least two}. Thus we may restate Question \ref{q:29} by asking whether a non-constant quasiregular mapping $\Sph^n \to \Sph^n$ either is a local homeomorphism or has branch set of topological dimension $n-2$. 

Also this question is connected to Questions \ref{q:12}--\ref{q:15}. Indeed, the double suspension $\Sigma^2 \varphi \colon \Sigma^2 \Sph^3 \to \Sigma^2 H^3$ of the covering map $\varphi \colon \Sph^3 \to H^3$ of the Poincar\'e homology sphere is a branched cover having one-dimensional branch set. To our knowledge, solutions to Questions \ref{q:27}, \ref{q:28}, and \ref{q:29} even for this particular map are not known. 

Regarding Question \ref{q:29}, we note in passing that, for $n=3$, the structure of the branch set of a (topological) branched cover $\Sph^3 \to \Sph^3$ is not understood. Especially, a classical question of Church and Hemmingsen \cite{Church-Hemmingsen} asks whether \emph{there exists a light interior map $f\colon M \to N$ between $3$-manifolds for which the topological dimension of $f(B_f)$ is zero}. 
For a light mapping in question, $\dim f(B_f) = \dim B_f$ by a result of Church and Hemmingsen in \cite{Church-Hemmingsen}. To our knowledge the question of Church and Hemmingsen is open also for branched covers.

Question \ref{q:30} has its origins in the classification of closed quasiregularly elliptic $2$ and $3$ manifolds and a question of Gromov \cite{Gromov-Hyperbolic}. The term \emph{quasiregularly elliptic $n$-manifold} was coined by Bonk and Heinonen in \cite{bh:cohomology} to denote a connected, oriented and Riemannian $n$-manifold $N$ which admits a non-constant quasiregular mapping $\R^n \to N$. Since all Riemannian metrics of a closed smooth manifold $N$ are quasiconformally equivalent, that is, $\id \colon (N,g_N) \to (N,g'_N)$ is quasiconformal for Riemannian metrics $g_N$ and $g'_N$, quasiregular ellipticity of a closed smooth $n$-manifold depends only on the topology (and possibly the smooth structure) of the manifold.

By the classical uniformization theorem for surfaces, the only closed quasiregularly elliptic $2$-manifolds are $\Sph^2$ and $\mathbb T^2$. By a theorem of Jormakka \cite{Jormakka}, the only closed quasiregularly elliptic $3$-manifolds are $\Sph^3$, $\Sph^2\times \Sph^1$, and $\mathbb T^3$ and their quotients; Jormakka's proof assumed validity of Thurston's Geometrization conjecture, later proved by Perelman.

In this terminology, Gromov asked in \cite[p.200]{Gromov-Hyperbolic} whether there exist simply connected closed manifolds that are not quasiregularly elliptic. At the time of Question \ref{q:30}, spheres were the only examples in the quasiregular literature of simply connected quasiregularly elliptic manifolds. Question \ref{q:30}, promoted by Rickman, was solved also by him in \cite{rick:qr-elliptic-4-manifolds} by an elaborate construction in which two restrictions of 
Alexander maps $\mathbb T^2 \times \mathbb T^2 \to \Sph^2 \times \Sph^2$ are glued together using a method of Piergallini \cite{Piergallini}.

\begin{theorem}[Rickman]
There exists a quasiregular mapping $\mathbb T^4 \to (\Sph^2\times \Sph^2)\# (\Sph^2\times \Sph^2)$. A fortiori, there exists a quasiregular mapping $\R^4 \to (\Sph^2\times \Sph^2)\# (\Sph^2\times \Sph^2)$.
\end{theorem}

Gromov's question was answered to the positive by Prywes \cite{pry:cohomology}: \emph{If $N$ is a quasiregularly elliptic $n$-manifold, then $\dim H_{\mathrm{dR}}^k(N) \le { n \choose k }$}, where $H^k_{\mathrm{dR}}(N)$ is the $k$th de Rham cohomology of $N$. Work of Heikkil\"a and Pankka improved the conclusion in Prywes' theorem to an algebraic result \cite{Heikkila-Pankka}: \emph{The de Rham cohomology $H_{\mathrm{dR}}^*(N)$ of a closed quasiregularly elliptic $n$-manifold admits an algebra monomorphism $H_{\mathrm{dR}}^*(N) \to \wedge^* \R^n$}. Together with a result of Piergallini and Zuddas \cite{Piergallini-Zuddas} on the existence of branched covers $\mathbb T^4 \to N$ into simply connected $4$-manifolds, and the classification of simply connected $4$-manifolds, this cohomology embedding result yields a  topological classification of closed simply connected quasiregularly elliptic $4$-manifolds. This classification was extended to a full classification of closed quasiregularly elliptic $4$-manifolds by Manin and Prywes \cite{Manin-Prywes}.

\begin{theorem}[Heikkil\"a--Pankka, Manin--Prywes]
A closed $4$-manifold $N$ is quasiregularly elliptic (with respect to at least one smooth structure) if and only if $N$ is homeomorphic to a quotient of $\#^k (\Sph^2\times \Sph^2)$, $\#^i \mathbb{C}P^2 \#^j \overline{\mathbb{C}P}^2$, or $\Sph^{4-k} \times (\Sph^1)^k$, where $0\le k,i,j\le 3$.
\end{theorem}

Regarding Question \ref{q:30}, we would also like to highlight a parallel development regarding elliptic manifolds. As coined by Gromov (see \cite[Section 2.41]{Gromov}, a manifold $M$ is \emph{elliptic} if there exists a Lipschitz mapping $\R^n \to M$ of positive asymptotic degree. Berdnikov, Guth, and Manin \cite{Berdnikov-Guth-Manin} showed that \emph{a formal, simply connected, and closed manifold is elliptic if and only if it is scalable}. We refer to \cite{Berdnikov-Guth-Manin} for a detailed discussion on terminology, and merely say here that a closed manifold $M$ is \emph{formal} if it admits a self-map $M \to M$ having the property that the induced map in cohomology $H^k_{\dR}(M) \to H^k_{\dR}(M)$ is given by multiplication with $p^k$, for some $p\in \N$. Coined by Manin and Berdnikov \cite{Berdnikov-Manin}, a formal closed manifold $M$ is \emph{scalable} if there exists a constant $C>0$ and a logarithmically dense set of $p\in \N$ for which there exists a $C(p+1)$-Lipschitz self-map which induces a multiplication by $p^k$ on cohomology $H^k_{\dR}(M)$. 

The parallels between elliptic and quasiregularly elliptic theory are highlighted in a recent result of Berdnikov, Guth, and Manin \cite[Theorem C]{Berdnikov-Guth-Manin} stating that a formal, simply connected, and closed $n$-manifold $M$ is elliptic if and only if there is an embedding of algebras $H^*_{\dR}(M) \to \wedge^* \R^n$. Note that, contrary to the quasiregularly elliptic case, the existence of the embedding of algebras is also sufficient. This follows from the construction of $1$-Lipschitz maps $\R^n \to M$ from topological data (see \cite[Theorem 5.1]{Berdnikov-Guth-Manin}). To our knowledge, similar constructions for quasihregular mappings are not known. Also the quasiregular side of the qualitative homotopy theory, developed in these papers, seems still to be understood.

\section{Geometric topology redux}\label{sec:geometric-topology-2}

\begin{question}\label{q:31}
\it Does every closed oriented topological $n$-manifold, $n \ge 4$, admit a branched cover onto $\Sph^n$?
\end{question}

This question remains open.

\begin{question}\label{q:32}
\it Does every closed oriented PL $n$-manifold, $n \ge 5$, admit a PL branched cover of degree at most $n$ onto the PL $n$-sphere?
\end{question}

This question remains open.

\begin{question}\label{q:33}
\it Does every closed topological four-manifold admit a metric that makes the manifold Ahlfors $4$-regular and linearly locally contractible?
\end{question}

This question remains open.

\subsection{Definitions and notation}\label{subsec:geometric-topology-2-definitions-and-notation}

A piecewise linear (PL) structure on a manifold $M$ is an atlas of charts whose transition mappings are piecewise linear. A manifold $M$ with a piecewise linear structure is called a \emph{PL} manifold; we refer to the monograph of Rourke and Sanderson \cite{Rourke-Sanderson} for a detailed discussion on PL topology. We also recall that an atlas with quasiconformal transition mappings is called a \emph{quasiconformal structure} and a manifold with a quasiconformal structure is called a \emph{quasiconformal manifold}. 

\subsection{Commentary}\label{subsec:geometric-topology-2-commentary}

These questions, ending the list of Heinonen and Semmes, are inherently questions on structures, given by atlases, on manifolds. To our knowledge, there has been no progress on these questions and we have very little to add to the discussion of Heinonen and Semmes on these problems.

In low dimensions $n\le 3$, each manifold admits a PL structure; for $1$-manifolds this is elementary, for $n=2$ this follows from the classification of surfaces and is due to Rad\'o and Poincar\'e, and for $n=3$ this follows from Moise's {\it Hauptvermutung} \cite{Moise}. For $n=4$, there are closed orientable topological manifolds that do not admit a quasiconformal structure by a theorem of Donaldson and Sullivan \cite{Donaldson-Sullivan}; in fact, Donaldson and Sullivan show that a closed simply connected $4$-manifold with intersection form $-E_8$ does not have a quasiconformal structure. For $n\ge 5$ each topological manifold admits a unique quasiconformal -- even Lipschitz -- structure by a theorem of Sullivan \cite{Sullivan}.

Regarding Question \ref{q:31}, it is a classical theorem of Alexander \cite{Alexander} that each orientable PL $n$-manifold $M$ admits a PL branched cover $M\to \Sph^n$. Thus Question \ref{q:31} may be viewed as asking for a topological version of Alexander's theorem. To our knowledge, this question is open in all dimensions $n\ge 4$. It should be noted that, in dimensions $n\le 3$, there are robust methods to construct branched covers between manifolds with boundary (see Berstein and Edmonds \cite{ Berstein-Edmonds_TAMS, Berstein-Edmonds_IJM}) and methods to deform non-zero degree maps to branched covers (see Edmonds \cite{Edmonds_MA, Edmonds_Annals}). Heinonen and Rickman used these methods in \cite{Heinonen-Rickman_Duke} to construct a bounded quasiregular map $B^3 \to B^3$ without radial limits in a limit set of a Kleinian group. More precisely they prove the following theorem.

\begin{theorem}[Heinonen--Rickman {\cite[Theorem 9.3]{Heinonen-Rickman_Duke}}]
Let $\Gamma$ be a geometrically finite torsion free Kleinian group without parabolic elements acting on $\mathbb S^2$ with limit set $\Lambda_\Gamma$ not the whole sphere. Then there exists a bounded quasiregular mapping $f\colon B^3 \to B^3$ such that $f$ has no radial limit points at points $\Lambda_\Gamma$. 
\end{theorem}
Recall that classical Fatou's theorem gives that a holomorphic map $\mathbb D \to \mathbb D$ has radial limit almost everywhere. It remains unknown whether a bounded quasiregular map $B^3 \to B^3$ has a single radial limit; see, however, result of K.~Rajala \cite{Rajala_AJM} showing that for locally homeomorphic quasiregular mappings $B^n \to \R^n$ has radial limits.

Constructions for quasiregular mappings in dimensions $n=4$ are more sparse and we refer here to Piergallini \cite{Piergallini} and works of Piergallini and Zuddas \cite{Piergallini-Zuddas-JLMS, Piergallini-Zuddas} on constructions of piecewise linear branched covers. 

Regarding Question \ref{q:32}, the answer is affirmative in dimensions $n\le 4$. For $n=2$, since each closed orientable surface is a connected sum of tori by the classification of surfaces, each such surface admits an orientation preserving involution $\iota \colon \Sigma \to \Sigma$, having a discrete fixed point set, for which the corresponding quotient space $\Sigma/\iota$ is a $2$-sphere and the quotient map $\Sigma \to \Sigma/\iota$ is a branched cover of degree $2$. In dimension $n=3$, every $3$-manifold admits a PL structure by Moise's theorem \cite{Moise} and the Hilden--Hirsh--Montesinos theorem \cite{Hilden-BAMS,Hilden-AJM, Hirsh, Montesinos} states that \emph{each closed and oriented $3$-manifold $M$ admits a degree $3$ PL branched cover $M\to \Sph^3$}. In dimension $n=4$, by a theorem of Piergallini \cite{Piergallini-4fold}, \emph{each closed and orientable PL $4$-manifold $M$ admits a degree $4$ PL branched cover $M\to \Sph^4$}.

Finally, Question \ref{q:33} may be interpreted as asking whether, although some manifolds may not carry a geometrically meaningful atlas, perhaps all manifolds nonetheless carry a meaningful metric -- in spirit of Questions \ref{q:3}, \ref{q:4}, \ref{q:5}, \ref{q:10}, and \ref{q:11} -- which allows first order analysis; see e.g.~Heinonen \cite{Heinonen-BAMS} and Semmes \cite{sem:selecta}.

\section{Concluding remarks}\label{sec:conclusion}

As noted in the introduction, the field of analysis on metric spaces has grown in new and unexpected directions since the original formulation of the Heinonen--Semmes questions. Without giving an exhaustive survey, we mention some of those that --- in our opinion --- are connected to the problems discussed here, and that have not yet been discussed in the commentaries above.

The emerging subject of analysis in metric spaces benefited tremendously from the appearance of Heinonen's book \cite{hei:lams}. Published in 2001, this book was based on notes from a graduate course taught by Heinonen\footnote{and attended by the third author as an aspiring beginning graduate student.} in the fall of 1996, and thus has its origins around the same time as the Heinonen--Semmes problem list. In this book, Heinonen presents the foundations of the newly developed subject with a focus on the theory and structure of metric spaces supporting a Poincar\'e inequality and/or the Loewner condition, basic elements of quasisymmetric mapping theory, and results on embeddability and metric deformations such as Assouad's embedding theorem and the David--Semmes deformation theory.

One of the first major developments in the field that came \emph{after} the publication of the Heinonen--Semmes list was Cheeger's 1999 paper \cite{cheeger} generalizing Rademacher's theorem (a.e. differentiability of Lipschitz functions on $\mathbb{R}^n$) to all doubling metric measure spaces supporting a (weak) Poincar\'e inequality, in the sense described in Section \ref{sec:mms}. The fact that such a strong theorem could be proved under such general assumptions, soon after the surprising new examples of Question \ref{q:19} and \ref{q:20} were constructed, led to extensive work on the structure of spaces supporting Poincar\'e inequalities. Cheeger's results immediately imply certain general non-embedding theorems; for instance, spaces as in Question \ref{q:20} cannot be subsets of any Euclidean space, even after a bi-Lipschitz change of the metric. This line of research was furthered substantially in work of Cheeger and Kleiner; see \cite{cheegerkleiner:banach} and \cite{cheegerkleiner:banach2} for the case of embeddings to Banach spaces with the Radon--Nikodym property, and \cite{cheegerkleiner:L1} and \cite{cheegerkleiner:L1-2} for the case of embeddings into the Banach space $L^1$. The main result of \cite{cheegerkleiner:L1}, in conjunction with work of Lee and Naor \cite{ln:gl}, showcases the sub-Riemannian Heisenberg group as a geometrically natural counterexample to a well-known conjecture of Goemans and Linial in algorithmic computer science, concerning computability of algorithms related to the Sparsest Cut Conjecture. These developments opened a rich and ongoing line of research at the intersection between analysis in metric spaces, metric embedding theory, and computer science. We refer to Naor's article in the proceedings of the 2010 ICM \cite{naor:icm} for more information. Recent work in this area has drawn heavily on advances in the study of uniform rectifiability in sub-Riemannian spaces; we discuss this point later in this section.

Other approaches to differentiation in non-smooth metric spaces, with interesting connections to Poincar\'e inequalities and the aforementioned work of Cheeger, have also advanced in recent years. One such approach originates in work of Weaver \cite{weaver}, who proposes a theory of \emph{derivations} on metric measure spaces: linear operators acting on Lipschitz functions that satisfy the functional analytic properties of classical partial derivatives. Another is the theory of \emph{Alberti representations}: certain decompositions of measures along fragments of rectifiable curves, reminiscent of Fubini's theorem. Alberti representations were studied intensively by Bate \cite{bate}.\footnote{See also the related earlier notion of \emph{test plans} defined in \cite{ags}.} Among other successes, Bate's work yields a new proof of Cheeger's differentiability theorem. The circle of ideas connecting Cheeger differentiation, Weaver derivations, and Alberti representations has been studied by Gong \cite{gong}, Schioppa \cite{schioppa}, and many others. Interactions between these theories and Poincar\'e inequalities is also now an extensive story that still contains a number of open questions; see work of Bate--Li \cite{bateli}, Eriksson-Bique \cite{eriksson-bique}, and Schioppa \cite{schioppa_example}.

Along with differentiation theory, many other aspects of classical geometric measure theory have made their way to the metric setting. One other major example is the theory of \emph{currents}, essentially, generalized surfaces defined as linear functionals on spaces of differential forms. In Euclidean space, the study of these objects was initiated by de Rham and extended by Federer and Fleming, and they are a foundational tool in variational calculus. The work of Ambrosio--Kirchheim \cite{ambrosiokirchheim} and Lang \cite{lang} brought currents into the metric measure setting, with (a heavily simplified statement of) the main idea being to replace differential forms by tuples of Lipschitz functions. Schioppa \cite{schioppa_currents} relates these objects to the differentiation theories mentioned above. A recent beautiful application of metric currents, germane to Section \ref{sec:mms} of this manuscript, appears in the work of Basso--Marti--Wenger \cite{bassomartiwenger}, who use metric currents to (among other things) give a new proof of a result of Semmes \cite{sem:selecta} on the existence of Poincar\'e inequalities on metric manifolds.\footnote{Much more could be said about metric currents in their own right, independent of connections to differentiability, Poincar\'e inequalities, and the main topics of the Heinonen--Semmes list, but we forgo such a discussion here.}

The theory of first-order Sobolev spaces on metric measure spaces has its origins in work of Korevaar--Schoen \cite{ks:sobolev} and Haj{\l}asz \cite{haj:sobolev} predating the Heinonen--Semmes problem list and the appearance of \cite{hei:lams}. However, the appearance of the Newtonian--Sobolev space, which was introduced in Cheeger's differentiation paper \cite{cheeger} and also in the thesis of Shanmugalingam \cite{shan:sobolev}, brought a new perspective by grounding the theory of such function spaces in the notion of modulus of curve families, thereby drawing ties to existing and ongoing work in quasisymmetric mapping theory. The theory of the Newtonian--Sobolev spaces is nowadays a rich and fruitful area of independent interest, and in order to keep our disussion brief we do not go into details here. Much work in this area is due to Shanmugalingam, J. Bj\"orn, and A. Bj\"orn, and their students and collaborators. The interested reader is invited to consult the books \cite{bj:book} and \cite{hkst:book} for further information. We remark that the Haj{\l}asz--Sobolev space introduced in \cite{haj:sobolev} has also been the subject of intense study; this framework is particularly suited for generalizations to fractional order Sobolev spaces and retains nontrivial interest on fractal spaces which fail to support sufficiently rich families of rectifiable curves to make the Newtonian--Sobolev theory nontrivial. The relationship between these notions of Sobolev space, and others such as the spaces of Korevaar--Schoen, is explained in detail in \cite[Chapter 10]{hkst:book}.

Another important early success of the newly developed subject consisted of consequences for geometric group theory and dynamical systems, especially via the lens of {\it conformal dimension}. In particular, Bonk and Kleiner's solutions to questions \ref{q:3} and \ref{q:5} on the list showcased the application of analysis in metric spaces to rigidity and uniformization problems in the geometry of infinite groups. Informative overviews of these developments can be found in the articles by Bonk and by Kleiner, published in the proceedings of the 2006 ICM \cite{bonk:icm}, \cite{kleiner:icm}. 

These results are the tip of an iceberg of further questions concerning quasiconformal classification and uniformization of non-smooth sets and spaces. We remarked already in Section \ref{subsec:qc-uniformization-parametrization-commentary} that the analog of Question \ref{q:3} is known to have a positive answer for spheres in dimensions $1$ and $2$ but a negative answer for higher-dimensional spheres. Stepping away from the manifold case, one may ask about the quasiconformal geometry of classical fractals such as Sierpi\'nski carpets. Bonk's work \cite{bonk} shows that the complement of any collection $\{C_i\}$ of disjoint closed Jordan domains in the Riemann sphere is quasisymmetric to a `round carpet'' (i.e., a topological carpet whose complementary domains are round circles), provided only that the boundaries $\{\partial C_i\}$ satisfy a scale-invariant separation condition and are uniformly quasisymmetrically equivalent to $\mathbb{S}^1$. On the classification side, the work of Bonk--Merenkov \cite{bonkmerenkov, bonkmerenkov2} shows that all the classical self-similar `square Sierpi\'nski carpets' are quasisymmetrically distinct and have finite quasisymmetric symmetry groups. To the best of our knowledge, analogous conclusions for more complicated fractals like the (standard, cubical) Menger sponge remain wide open. (One basic stumbling block there is that the peripheral circles in a Sierpi\'nski carpet are topologically distinguished --- this is essential in \cite{bonkmerenkov} --- while by contrast the Menger sponge is topologically homogeneous in a very strong sense.) In this direction, however, we note that the exotic metric spaces of Bourdon--Pajot and Laakso, discussed in connection with Questions \ref{q:19} and \ref{q:20}, are topologically equivalent to the standard Menger sponge and enjoy strong self-similarity properties, but are not quasisymmetrically equivalent to the standard Menger sponge, since all of these spaces are minimal for conformal dimension and the Menger sponge is not.
The preceding circle of ideas is closely connected to developments in geometric group theory and complex dynamics. For instance, \cite{bonkmerenkov2} and \cite{merenkov:local} show that the classical square Sierpi\'nski carpet cannot be quasisymmetric either to the Julia set of a post-critically finite rational map or to the boundary of a torsion-free hyperbolic group. See also \cite{bonklyubichmerenkov}. 

Recent years have witnessed an explosion of interest in dynamics on nonsmooth metric spaces, especially nonsmooth metric structures on the $2$-sphere. This line of research was initiated in an extensive 
book by Bonk and Meyer \cite{Bonk-Meyer_book} and a monograph of Ha\"{\i}ssinsky and Pilgrim \cite{hp:book}. To briefly explain the main ideas, we recall that each expanding Thurston map $f\colon \mathbb S^2 \to \mathbb S^2$ induces a visual metric $\rho$ on $\mathbb S^2$ under which $f$ is a Lipschitz map. The construction is analogous to the construction of visual metrics on the boundary of a Gromov hyperbolic group, and also is allied to the notion of {\it hyperbolic filling} of a compact metric space. References for the latter topic include \cite{bp:besov} and \cite{bs:Sobolev-and-hyperbolic-filling}. Although the dynamically defined visual metrics $\rho$ of Bonk and Meyer need not be doubling, they provide an interesting class of metrics relevant for the quasisymmetric uniformization problem in Question \ref{q:3}. Indeed, in view of \cite[Theorem 18.1]{Bonk-Meyer_book} and \cite[Theorem 4.2.11]{hp:book}, \emph{the metric sphere $(\mathbb S^2, \rho)$ is quasisymmetric to the standard sphere $\mathbb S^2$ if and only if the expanding Thurston map $f$ is topologically conjugate to a rational map}.

Another development in dynamics associated to geometric mapping theory, dating back to the time of the Heinonen--Semmes problem list, is the introduction of uniformly quasiregular mappings by Iwaniec and Martin on Riemannian $n$-manifolds in \cite{Iwaniec-Martin}; see also Sullivan \cite{Sullivan-1983} for uniform quasiregularity in two dimensions. The rigidity of higher dimensional conformal geometry, as demonstrated by the classical Liouville's theorem, gives the impression that higher dimensional (quasi)conformal dynamics reduces to the theory of Kleinian groups. In their seminal work, Iwaniec and Martin showed that there are non-injective quasiregular mappings $\mathbb S^n \to \mathbb S^n$ whose distortion stays bounded under iteration. Such maps are called \emph{uniformly quasiregular (UQR)}. Martin and Peltonen \cite{Martin-Peltonen} showed that these mappings in fact exist in abundance on the $n$-sphere.
Constructions of Mayer \cite{Mayer_CDG}, Peltonen \cite{Peltonen}, Astola--Kangaslampi--Peltonen \cite{Astola-Kangaslampi-Peltonen}, Kangaslampi--Peltonen--Wu \cite{Kangaslampi-Peltonen-Wu}, and Fletcher--Wu \cite{Fletcher-Wu} show that spaces covered by a sphere admit a rich UQR mapping theory. On the other hand, works of Kangasniemi \cite{Kangasniemi_AJM, Kangasniemi_Adv} show that the closed manifolds admitting non-branching quasiregular dynamics form a strict subclass of quasiregularly elliptic manifolds. It is interesting to observe how the quasiregular results of Kangasniemi bear similarity to aforementioned works of Berdnikov, Guth, and Manin \cite{Berdnikov-Guth-Manin} in the Lipschitz category.

The late 1990's and early 2000's saw a surge of results relaxing the boundedness of the distortion in higher dimensional mapping theory. The two-dimensional result of David \cite{David-1988} on the Measurable Riemann Mapping Theorem for mappings whose Beltrami tensor may reach up to the unit circle was an antecedent to the study of \emph{mappings of finite distortion}. Regarding this vast field with connections to nonlinear elasticity, we refer here 
refer to the monographs of Iwaniec--Martin \cite{Iwaniec-Martin-book-2001} and Hencl--Koskela \cite{Hencl-Koskela-2014} for comprehensive presentations of the theory. Very recently, Kangasniemi and Onninen have found a completely new point of view on higher dimensional geometric mapping theory involving a new notion of \emph{quasiregular values}. In this direction, Kangasniemi and Onninen have obtained a version of the aforementioned Reshetnyak's theorem \cite{Kangasniemi-Onninen-TAMS} and Rickman's Picard theorem \cite{Kangasniemi-Onninen-GAFA}, both without pointwise quasiconformality assumptions on the mappings.

We next turn to the interplay between analysis in metric spaces and sub-Riemannian geometry, especially in the theory of sub-Riemannian metrics on the Heisenberg group and other nilpotent stratified Lie groups (Carnot groups). The ongoing development of these two subjects over the past forty years have been closely linked: sub-Riemannian spaces, equipped with their canonical {\it Carnot--Carath\'eodory metrics} form a natural testing ground for questions in analysis in metric spaces, while the additional differentiable and algebraic structure inherent in Carnot groups allows for a much more detailed study along the lines of Euclidean and Riemannian geometry. 

The field of sub-Riemannian geometry has a long history, eventually weaving together disparate lines of research in radically different fields including non-holonomic mechanics, several complex variables, quantum mechanics, nilpotent group theory, and geometric rigidity. From the viewpoint of its eventual alignment with analysis in metric spaces, Mostow's signature results on the rigidity of hyperbolic manifolds \cite{mos:rigidity1}, \cite{mos:rigidity2}, \cite{mos:rigidity3} are of special importance. Mostow's remarkable work, and in particular the use of quasiconformal mapping theory on the boundaries of the universal covers, provided substantial motivation for the development of quasiconformal mapping theory in the Heisenberg group and other Carnot groups arising as local models for the Gromov boundaries of rank one symmetric spaces of noncompact type. This line of research played an essential role in the eventual rise of analysis in metric spaces as a subject; we refer the reader to \cite{hei:calculus-in-carnot} and to the introduction of \cite{kr:foundations} for a pr\'ecis of this compelling story. 

The problems in the Heinonen--Semmes list presented in section \ref{sec:heisenberg} bear the stamp of these historical precedents. For instance, Questions \ref{q:21} and \ref{q:25} (which remain open) concern quasiconformal and quasisymmetric mapping theory in the sub-Riemannian Heisenberg setting; the latter problem addresses the potential flexibility of quasiconformal self-maps of the Heisenberg group, while the former asks to what extent the metric doubling deformation theory developed by David and Semmes --- which in particular gives rise to quasisymmetric mappings --- transfers to the Heisenberg group. 

Numerous other lines of research in sub-Riemannian analysis and geometry were already taking shape at the time of the Heinonen--Semmes survey, and the field has expanded tremendously since that time. We can do no more here than point the interested reader to some helpful references for further information. With regards to rectifiability in sub-Riemannian spaces, the foundational work of Ambrosio--Kirchheim \cite{ambrosiokirchheim} indicated the non-applicability of classical notions of rectifiability outside of sufficiently low-dimensional objects. A novel approach to intrinsic sub-Riemannian rectifiability was begun by Franchi--Serapioni--Serra Cassano \cite{fss}, see also Serra Cassano's survey \cite{serracassano} for an informative introduction and overview. A related line of research focuses on the intrinsic differential geometry of classically smooth submanifolds, drawing motivation from sub-Riemannian analogs of well-known questions in the calculus of variations. A signature open question in this direction is Pansu's isoperimetric conjecture \cite{pan:isoperimetric1}, \cite{pan:isoperimetric2} on the classification of minimizers of the sub-Riemannian perimeter functional under a volume constraint. We direct the reader to \cite{cdpt:book} for an outline of the state of the art with regards to partial results towards Pansu's conjecture as of 2007, but note that developments since 2007 have advanced that field further.

Many of the previous directions come together in recent advances in the theory of uniform rectifiability in sub-Riemannian and metric spaces. As mentioned in Section \ref{sec:bl-parametrizations-embeddings}, the theory of uniform rectifiability initiated by David and Semmes in \cite{ds:uniform-rectifiability} is by now a major program in non-smooth analysis. Questions \ref{q:9}, \ref{q:22}, \ref{q:23}, and \ref{q:24} above were early attempts to advance this sort of quantitative geometric measure theory in non-Euclidean spaces. The outcomes have been especially fruitful in Carnot groups, and particularly in the Heisenberg group which has seen some spectacular results and unexpected connections to the aforementioned applications to theoretical computer science. We direct the reader to \cite{cheegerkleinernaor}, \cite{naoryoung1}, and \cite{naoryoung2} for more information.

We conclude this historical overview with brief remarks on other areas of current interest, beginning with synthetic differential geometry. This field has a long history dating back to the work of Alexandrov, Gromov, and others who captured the intrinsic metric content of bounds on sectional curvature. By the time that the Heinonen--Semmes list appeared, it was already recognized that Gromov--Hausdorff limits of Riemannian manifolds with suitable Ricci curvature bounds could supply other examples of doubling spaces supporting Poincar\'e inequalities \cite{cheegercolding}. Much later, Lott--Villani \cite{lottvillani} and Sturm \cite{sturm-acta-I}, \cite{sturm-acta-II} used optimal transport to give `synthetic' definitions of lower Ricci bounds in metric measure spaces, from which Poincar\'e inequalities can also be deduced \cite{tapiorajala}. An interesting byproduct of this line of study was an intrinsic proof for the reflexivity of the Newtonian--Sobolev spaces avoiding the use of the Poincar\'e inequality \cite{acdm:reflexivity}. This topic has also seen intense progress in the last few decades. The curious reader is invited to consult \cite{ags:book}, \cite{ags}, \cite{gig:memoirs1}, \cite{gig:memoirs2}, or \cite{gp:nonsmooth} for more information on this tremendously active field.

On the other hand, a field that was for many years largely separate from `analysis on metric spaces' as understood here was the subject of `analysis on fractals'. This area is concerned with topics such as Dirichlet energy and heat kernel estimates on fractal sets. We refer the reader to \cite{strichartz}, \cite{kigami}, and {\cite{barlow}. Recent years have seen the rise of novel connections between the two fields, which are now much more closely allied. For example, we highlight Murugan's discovery \cite{murugan} that Hino's martingale construction \cite{hino} recovers the dimension of the Cheeger Lipschitz differentiability structure.

Finally, another topic that is mostly absent from the problem list and the preceding discussion is probability. More work has been done recently on \emph{random geometry}: the structure of (typically non-smooth) random geometric objects. Well-known examples include (traces or graphs of) Brownian motion and the Schramm--Loewner evolution, but there are also natural random metrics on the sphere (the ``Brownian map'') \cite{legall} and on trees \cite{aldous}. In the vein of Section \ref{sec:qc-uniformization-parametrization}, there are many interesting questions to be asked about quasiconformal uniformization and distortion of such random structures. For samples, see \cite{troscheit}, \cite{binderhakobyanli}, and \cite[Problem 9.1]{bonkmeyer}. We anticipate that the continued development of stochastic analogs of many of the preceding lines of research will keep researchers busy for many years to come.

Despite the length of this final section, we have only scratched the surface in our description of the rich literature comprising the modern theory of analysis in metric spaces. While the subject owes a debt to many historical precursors, the problems posed by Heinonen and Semmes provided generations of new researchers with a diverse and compelling array of questions. Efforts towards answing these questions spurred the development of many signature results in the field, and helped to bring the subject to where it is today. We hope that this survey article has provided our readers with a new perspective on a field which continues to flourish and thrive.

\bibliographystyle{acm}

\bibliography{biblio}
\end{document}